\documentstyle[11pt]{article}
\begin{document}
\author{{Shiqiu Zheng\footnote{E-mail: shiqiumath@163.com.}}\\
  \\
\small College of Sciences, Tianjin University of Science and Technology, Tianjin 300457, China
}
\date{}
\title{\textbf{On $g$-expectations and filtration-consistent nonlinear expectations}\thanks{This work is supported by the Foundation of Ministry of Education of China (No. 20YJCZH245).}}\maketitle

\textbf{Abstract:}\ \ In this paper, we obtain a comparison theorem and a invariant representation theorem for backward stochastic differential equations (BSDEs) without any assumption on the second variable $z$. Using the two results, we further develop the theory of $g$-expectations. Filtration-consistent nonlinear expectation (${\cal{F}}$-expectation) provides an ideal characterization for the dynamical risk measures, asset pricing and utilities. We propose two
new conditions: an absolutely continuous condition and a (locally Lipschitz) domination condition. Under the two conditions respectively, we prove that any ${\cal{F}}$-expectation can be represented as a $g$-expectation. Our results contain a representation theorem for $n$-dimensional ${\cal{F}}$-expectations in the Lipschitz case, and two representation theorems for $1$-dimensional ${\cal{F}}$-expectations in the locally Lipschitz case, which contain quadratic ${\cal{F}}$-expectations.    \\


\section{Introduction}
As typical time-consistent nonlinear expectations, $g$-expectations were initiated in Peng \cite{Peng97} via the solutions of backward stochastic differential equations (BSDEs) in the Lipschitz case. The theory of $g$-expectations has been widely applied in asset pricing, utility theory and risk measures (see \cite{R}, \cite{Jiang}, \cite{DP}, \cite{Xu} and the references therein). To deal with more general cases, the notion of $g$-expectations was generalized by Jia \cite{Jia} in the uniformly continuous case, and by Ma and Yao \cite{MY} in the quadratic growth case. One of objectives of this paper is to further develop the theory of $g$-expectations. To this end, we need to study the comparison theorem and the invariant representation theorem for BSDEs in more general situations.

As one of the most important properties of $1$-dimensional BSDEs, comparison theorem plays an important role in the theory of BSDEs and its applications. We refer to El Karoui et al. \cite{EPQ} for the comparison theorem for Lipschitz continuous generators, and to Kobylanski \cite{K} and Briand and Hu \cite{BH} for the comparison theorems for quadratic growth generators. Different from the $1$-dimensional case, the comparison theorem for $n$-dimensional BSDEs needs some special restrictive conditions on the structure of the generator $g$ (see Hu and Peng \cite{HP}). It is clear that comparison theorem implies the uniqueness of solution. In this paper, we show that the converse is also true. Specifically, we get a general comparison theorem for BSDEs without any assumption on the second variable $z$ (see Theorem 2.7), whose proof mainly relies on the uniqueness of solution and the comparison theorem of a special SDE. Our result means that comparison theorem usually holds true as long as the BSDE has a unique (resp. maximal) solution in some space.

Another important property of BSDEs is the invariant representation theorem for generators, which is a powerful tool to study the generators from the solutions of BSDEs. In order to study the converse comparison problem for BSDEs, the invariant representation theorem was firstly proved in Briand et al. \cite{BP}, and then generalized by \cite{Jiang05}, \cite{Jiang}, \cite{Jia}, \cite{FJ} and \cite{Zheng2}, etc. To our best knowledge, all of these studies rely on some proper estimates on the solutions of BSDEs, which need some assumptions on the second variable $z$. In this paper, we give an explicit solution to a special $n$-dimensional mean-field BSDE, from which a $n$-dimensional invariant representation theorem is derived for deterministic generators without any assumption on $z$ (see Proposition 2.8). Using this result, we obtain a converse comparison theorem for $n$-dimensional BSDEs. Interestingly, we also give the explicit solutions to two special $n$-dimensional mean-field BSDEs with $L_n^2$ terminal variables, without any assumption on $z$ (see Example 2.9).

The axiomatic system of filtration-consistent nonlinear expectation (${\cal{F}}$-expectation) was initiated in Coquet et al. \cite{CH}. We refer to Peng \cite{Peng04} for an equivalent definition. It is well-known that ${\cal{F}}$-expectation gives an ideal characterization for the dynamical behaviors of risk measure, asset pricing and utilities, and also can be considered as a dynamic version of the famous risk measures introduced by Artzner et al. \cite{Ar} and F\"ollmer and Schied \cite{FS} (see Jiang \cite{Jiang} for a $g$-expectations case). In \cite{CH}, the authors raised the problem: \emph{is this notion of $g$-expectations general enough to represent all "enough regular" ${\cal{F}}$-expectations?} This question was considered to be theoretically very interesting and practically important by Peng in his ICM lecture (see Peng \cite[Page 403]{Peng10}). The first result on this representation problem was obtained by \cite{CH}, which proved that any ${\cal{F}}$-expectation ${\cal{E}}$ is a $g$-expectation whose generator $g$ is independent of $y$ and Lipschitz in $z$,
when it satisfies the following domination condition
$${\cal{E}}[\xi|{\cal{F}}_t]-{\cal{E}}[\eta|{\cal{F}}_t]\leq{\cal{E}}^{\mu}[\xi-\eta|{\cal{F}}_t],\ \ \forall\xi,\eta\in L^2({\cal{F}}_T),\ \forall t\in[0,T],\ \eqno(1.1)$$
where ${\cal{E}}^\mu$ is the $g$-expectation defined via the BSDE with generator $g=\mu|z|$. Note that (1.1) is equivalent to the original domination condition plus the strict monotonicity condition and translation invariance condition in \cite{CH} (see Section 4). We refer to Cohen \cite{Co} and Royer \cite{Ro} for the corresponding results under (1.1) in general filtrations. Furthermore, under several domination conditions, Hu et al. \cite{Yao} proved that any ${\cal{F}}$-expectation with independent increments is a $g$-expectation whose generator $g$ is determinstic, independent of $y$, and has quadratic growth (locally Lipschitz) in $z$. Under a general domination condition ($\mu|z|$ in (1.1) is replaced by an increasing and linear growth function $\phi(|z|)$), Zheng and Li \cite{Zheng1} proved that any ${\cal{F}}$-expectation is a $g$-expectation whose generator $g$ is independent of $y$ and uniformly continuous in $z$. We point out that the studies in \cite{Ro}, \cite{Yao}, \cite{Co} and \cite{Zheng1} all use the domination method developed in \cite{CH}, which uses the domination conditions to obtain the convergence and the Doob-Mayer decomposition for nonlinear supermartingales. However, it is known from \cite{Yao} and \cite{Zheng1} that studying such problem by the domination method is by no means easy without (1.1). One reason is that it is not easy to find a domination condition for ${\cal{F}}$-expectations. In fact, the domination condition may not exist in some non-Lipschitz cases. Another reason is that even if the domination condition can be found, some new techniques and methods will have to be developed in the proof, since the $L^2$-estimates and $L^2$-dominations crucial in the proof of the representation theorem in \cite{CH} are usually not true in the non-Lipschitz cases.

Another objective of this paper is to study the above representation problem under independent increments condition. To this end, we propose two
new conditions: an absolutely continuous condition and a (locally Lipschitz) domination condition. Specifically, we answer this problem in two directions. The first one is that when the ${\cal{F}}$-expectation ${\cal{E}}$ satisfies the following absolutely continuous condition
$${\cal{E}}[zB_t]\ \textmd{is absolutely continuous in}\ t,\ \ \ \forall z\in{\bf{R}}^{n\times d},\eqno(1.2)$$
and certain convergence, we prove that any ${\cal{F}}$-expectation ${\cal{E}}$ is a $g$-expectation whose generator $g$ is deterministic, independent of $y$ and Lipschitz continuous in $z$ or locally Lipschitz continuous in $z$ (see Theorem 4.3). Our proof does not use the domination method mentioned above. Our result contains a representation theorem for $n$-dimensional ${\cal{F}}$-expectation in the Lipschitz case (Note that the representation theorems in \cite{CH}, \cite{Ro}, \cite{Yao}, \cite{Co} and \cite{Zheng1} are all for $1$-dimensional ${\cal{F}}$-expectations). This $n$-dimensional representation can not be fully obtained by the domination method, since the domination condition relies on the comparison theorem of BSDEs, but as mentioned above, in general, the $n$-dimensional comparison theorem are not true (see \cite{HP}). Our result also contains a representation theorem for quadratic ${\cal{F}}$-expectations with unbounded terminal variables, a notable example of which is the entropic risk measure (see Barrieu and El Karoui \cite{BEl} for a bounded terminal variables case). Moreover, we also show that (1.2) is a necessary condition for this representation problem (see Proposition 4.4). In our opinion, (1.2) is more easy to be checked in applications than (1.1) and the domination conditions in \cite{Yao} and \cite{Zheng1}. Our method begins with the special terminal variable set ${\cal{R}}$ (the set of all random variables: $y+z(B_v-B_u)$), and then consider more general terminal variable sets.

The second one is that when the ${\cal{F}}$-expectation ${\cal{E}}$ satisfies the following (locally Lipschitz) domination condition
$${\cal{E}}[\xi|{\cal{F}}_t]-{\cal{E}}[\eta|{\cal{F}}_t]\leq{\cal{E}}^{\rho(k)}[\xi-\eta|{\cal{F}}_t],\ \ \forall k>0,\ \forall\xi,\eta\in{\cal{R}}^k,\ \forall t\in[0,T],\eqno(1.3)$$
where $\rho(\cdot):{\mathbf{R_+}}\rightarrow{\mathbf{R_+}}$ is a nondecreasing function, we prove that any ${\cal{F}}$-expectation ${\cal{E}}$ is a $g$-expectation whose generator $g$ is deterministic, independent of $y$ and locally Lipschitz continuous in $z$ (see Theorem 4.6). (1.3) can be considered as a generalization of (1.1), and can be used to dominate the solutions of BSDEs with locally Lipschitz continuous generators. Our proof follows the spirit of \cite{CH}, but we uses a different strategy, since some crucial results under (1.1) are not true under (1.3). We also begins with the set ${\cal{R}}$ as in the first direction, which indeed determine the basic part of $g$-expectations. We establish a special Doob-Mayer decomposition for ${\cal{E}}$-supermartingale $\psi(t)+zB_t$ (see Lemma A.6), which is enough in our proof. Our method relies on the boundeness of solution $Z_t$ heavily, which is inspired by the works of Briand and Elie \cite{BE} and Cheridito and Nam \cite{CN} depending on Malliavin calculus. Since we uses the domination condition on ${\cal{R}}$ and consider a very special ${\cal{E}}$-supermartingale, some subtle techniques of stopping times used in \cite{CH}, \cite{Yao} and \cite{Zheng1} are omitted in our proof. Our result contains a representation theorem for quadratic ${\cal{F}}$-expectations, which is similar to the result obtained by \cite{Yao} under several domination conditions different from (1.3). Comparing with the domination conditions in \cite{Yao} (see Definition 3.8(1)-(3) and (H4) in \cite{Yao}), (1.3) seems to be more simple and not hard to be verified in applications (see Remark 4.8). We also provides a simple proof for the representation theorem in \cite{CH} in the deterministic generator case (see Corollary 4.7).

The paper is organized as follows. In Section 2, we will prove a general comparison theorem and a invariant representation theorem for BSDEs. In Section 3, we will further develop the $g$-expectations theory. In Section 4, we will study the representation problem for ${\cal{F}}$-expectations. In the appendix, we will give a special nonlinear Doob-Mayer decomposition.
\section{Some results on BSDEs}
Let ${{(B_t)}_{t\geq0}}$ be a $d$-dimensional standard Brownian motion defined a given complete probability space $(\Omega ,\cal{F},\mathit{P})$. Let $({\cal{F}}_t)_{t\geq 0}$ denote the natural filtration generated by ${{(B_t)}_{t\geq 0}}$, augmented
by the $\mathit{P}$-null sets of ${\cal{F}}$. Let the nondecreasing function $\rho(\cdot):[0,\infty)\longmapsto[0,\infty)$, and the positive constants $T$, $\mu$, and $\nu$ be given. Let ${\cal{P}}$ denote the $({\cal{F}}_t)_{t\geq 0}$-progressive sigma-field on $[0,T]\times\Omega$. Let ${\cal{T}}_{t,T}$ be the set of all the $({\cal{F}}_t)_{t\geq0}$-stopping times whose values belong to $[t,T]$. Let $|\cdot|$ denote the
Euclidean norm on $\textbf{R}^{n\times d}$ and $\|\cdot\|_\infty$ denote the essential supremum of random variables or stochastic processes. We define ${\mathbf{B}}_n(k):=\{z\in\textbf{R}^{n\times d}:|z|\leq k\},$ $k>0.$ Note that in this paper, all the equalities and inequalities for random variables hold true in the almost sure sense, and we denote $x=y$ (resp. $\geq, \leq$), if $x=(x^1,\cdots,x^n)^{\mathsf{T}},y=(y^1,\cdots,y^n)\in\textbf{R}^{n}$ satisfying $x^i=y^i$ (resp. $\geq, \leq$) for each $1\leq i\leq n$. We consider the following space of generators:

$\textbf{G}_n:=\{g(\omega,t,y,z):\Omega\times[0,T]\times{\mathbf{R}}^n\times{\mathbf{R}}^{n\times d}\longmapsto{\mathbf{R}}^n,$ such that $g$ is measurable with respect to ${\cal{P}}\otimes{\cal{B}}({\mathbf{R}}^n)\otimes{\cal{B}}({\mathbf{R}}^{n\times d})\};$

${\cal{G}}_n:=\{g\in\textbf{G}_n:$ $\forall y_i\in{\mathbf{R}}^n, i=1,2,$ $|{g}(t,y_1,\cdot)-{g}(t,y_2,\cdot)|\leq\mu|y_1-y_2|,\ dt\times dP$-$a.e.$ and $\forall z\in {\mathbf{R}}^{n\times d},\ \int_0^T|g(t,0,z)|dt<\infty\};$

${\cal{G}}_n^0:=\{g\in{\cal{G}}_n:$ $\forall y\in {\mathbf{R}}^n, g(t,y,0)=0$, $dt\times dP$-$a.e.\};$

${\cal{G}}_n^z:=\{g\in{\cal{G}}_n:$ $g$ is independent of $y$, $dt\times dP$-$a.e.\};$

$\widehat{{\cal{G}}}_n:=\{g\in{\cal{G}}_n:$ $g$ is deterministic, $dt$-$a.e.\};$

${\cal{L}}_n^{\rho}:=\{g\in{\cal{G}}_n:$ $\forall z_i\in {\mathbf{R}}^{n\times d}, i=1,2, |{g}(t,\cdot,z_{1})-{g}(t,\cdot,z_{2})|\leq\rho(|z_1|\vee|z_2|)|z_{\mathrm1}-z_{2}|, dt\times dP$-$a.e\}$;

${\cal{L}}_n^{\nu}:=\{g\in{\cal{G}}_n:$ $\forall z_i\in {\mathbf{R}}^{n\times d}, i=1,2, |{g}(t,\cdot,z_{1})-{g}(t,\cdot,z_{2})|\leq\nu|z_{\mathrm1}-z_{2}|, dt\times dP$-$a.e\}$;

We also use the following spaces of random variables and stochastic processes:

$L_n^0({\mathcal{F}}_t):=\{\xi:$ ${\mathbf{R}}^n$-valued, ${\mathcal{F}}_t$-measurable random variable$\};$

$L_n^p({\mathcal{F}}_t):=\{\xi:$ random variable in $L_n^0({\mathcal{F}}_t)$, such that ${{E}}\left[|\xi|^p\right]<\infty\},\ p\geq1;$

$L_n^{\infty}({\mathcal{F}}_t):=\{\xi:$ random variable in $L_n^0({\mathcal{F}}_t)$, such that $\|\xi\|_{\infty}<\infty\};$

${\mathcal{S}}_n:=\{(\psi_t)_{t\in[0,T]}:$ $\textbf{R}^n$-valued, continuous $({\cal{F}}_t)_{t\geq0}$-adapted
process$\}$;

${\mathcal{S}}^p_n:=\{(\psi_t)_{t\in[0,T]}:$ process in ${\mathcal{S}}_n$ such that ${{E}}\left
[{\mathrm{sup}}_{0\leq t\leq T}|\psi _t|^p\right]<\infty \},\ p\geq1;$

${H}^p_{n\times d}:=\{(\psi_t)_{t\in[0,T]}:$ $\textbf{R}^{n\times d}$-valued, $({\cal{F}}_t)_{t\geq 0}$-progressively measurable process such that $\int_0^T|\psi_t|^pdt
<\infty\},\ p\geq1;$

${\cal{H}}_{n\times d}^p:=\{(\psi_t)_{t\in[0,T]}:$ process in ${H}^p_{n\times d}$ such that ${{E}}\left[\left(\int_0^T|\psi_t|^2dt\right)^{\frac{p}{2}}\right]
<\infty \},\ p\geq1;$

${\cal{H}}^{BMO}_{n\times d}:=\{(\psi_t)_{t\in[0,T]}:$ process in ${H}^2_{n\times d}$ such that $\sup_{\tau\in{\cal{T}}_{0,T}}\left\|E[\int_\tau^T|\psi_t|^2dt|{\cal{F}}_\tau]\right\|_\infty<\infty \};$

${\cal{H}}_\rho^{BMO}:=\{(\psi_t)_{t\in[0,T]}:$ process in ${\cal{H}}^{BMO}_{1\times d}$ such that $\sup_{\tau\in{\cal{T}}_{0,T}}\left\|E[\int_\tau^T\rho(|\psi_t|)^2dt|{\cal{F}}_\tau]\right\|_\infty<\infty\};$

${\mathcal{H}}^{\infty}:=\bigcup_{k>0}{\mathcal{H}}^{\infty,k}$, where ${\mathcal{H}}^{\infty,k}:=\{(\psi_t)_{t\in[0,T]}:$ process in ${H}^1_{1\times d}$ such that $\|\psi\|_\infty\leq k\}, k>0;$

$l_{n\times d}^2(0,T):=\{(\psi_t)_{t\in[0,T]}:$ process in ${H}^2_{n\times d}$ is deterministic$\};$

$l_{n\times d}^\infty(0,T):=\{(\psi_t)_{t\in[0,T]}:$ process in $l^2_{n\times d}$ is bounded$\}.$

For a smooth random variable $\xi\in L_1^2({\mathcal {F}}_T)$, let $D_t\xi$ be the Malliavin derivative of $\xi.$ Let ${\mathcal{D}}^{1,2}$ be the closure of the set of all Malliavin differentiable random variables $\xi\in L_1^2({\mathcal {F}}_T)$ with respect to the norm $\|\xi\|_{1,2}=[E|\xi|^2+E\int_0^T|D_t\xi|^2dt]^\frac{1}{2}$ (see \cite[Section 1.2]{Na}). Moreover, we define

${\cal{R}}_n:=\bigcup_{k>0}{\cal{R}}_n^k$, where ${\cal{R}}_n^k:=\{y+z(B_v-B_u): u,v\in[0,T]\ \textmd{and}\ (y,z)\in\textbf{R}^n\times\textbf{B}_n(k)\}, k>0;$

$\overline{{\cal{R}}}:=\bigcup_{k>0}\overline{{\cal{R}}}^k$, where $\overline{{\cal{R}}}^k:=\{\xi:\xi\in{\mathcal{D}}^{1,2}$ such that $\sup_{t\in[0,T]}\|D_t\xi\|_{\infty}\leq k\},\ k>0.$\\

For $g\in\textbf{G}_n$ and $\xi\in L_n^0({\mathcal{F}}_T)$, we consider the following BSDE$(g,T,\xi)$:
$$Y_{t}=\xi+\int_{t}^T g(s,Y_s,Z_s)
ds-\int_t^T Z_sdB_s,\ \ \,\ \ t\in[0,T].\eqno(2.1)$$
\textbf{Definition 2.1} For $g\in\textbf{G}_n$ and $\xi\in L_n^0({\mathcal{F}}_T)$, the solution of the BSDE$(g,T,\xi)$ is a pair of processes $(Y_t,Z_t)\in{\mathcal{S}}_n\times H_{n\times d}^2$ such that $\int_{0}^T|g(t,Y_t,Z_t)|dt<\infty,$ and (2.1) holds true. \\

From Pardoux and R\u{a}scanu \cite[Proposition 3.12 and Theorem 3.17]{PR}, we can get\\\\
\textbf{Lemma 2.2} \emph{Let $g\in{\cal{G}}_n$ and $(y,Z_s)\in {\mathbf{R}}^n\times H_{n\times d}^2$ such that $\int_{0}^T|g(s,0,Z_s)|ds<\infty.$ Then the following SDE
$$Y^{y,Z,g}_t=y-\int_{0}^tg(s,Y_s^{y,Z,g},Z_s)
ds+\int_0^t Z_sdB_s,\ \ t\in[0,T],\eqno(2.2)$$
has a unique solution $Y^{y,Z,g}_t\in{\mathcal{S}}_n$. Moreover, if $n=1$, then for $\tilde{y}\in{\mathbf{R}}$ and $f_t\in H_{1\times1}^1$ such that $y\leq \tilde{y}$ and $g(t,Y^{\tilde{y},Z,f}_t,Z_t)\geq f_t$, $dt\times dP$-a.e., where $Y_t^{\tilde{y},Z,f}=\tilde{y}-\int_0^tf_sds+\int_0^tZ_sdB_s,\ t\in[0,T],$ we have for each $t\in[0,T]$, $Y^{y,Z,g}_t\leq Y^{\tilde{y},Z,f}_t$.}\\

Let ${\cal{Z}}_{n\times d}$ denote a subset of $H_{n\times d}^2$ such that for each $\eta_t,\zeta_t\in{\cal{Z}}_{n\times d}$ and $\tau\in{\cal{T}}_{0,T}$, we have $$\eta_t1_{[0,\tau)}(t)+\zeta_t1_{[\tau,T]}(t)\in{\cal{Z}}_{n\times d}.$$
For $g\in{\cal{G}}_n$, we define the following set:
$$\widetilde{{\cal{R}}}^g({\cal{Z}}_{n\times d}):=\{Y_\cdot^{y,Z,g}:\ y\in\textbf{R}^n\ \textmd{and}\ Z_t\in{\cal{Z}}_{n\times d}\ \textmd{such that}\ \int_0^T|g(t,0,Z_t)|dt<\infty\},$$
where $Y_\cdot^{y,Z,g}$ is the process given by (2.2). We set $\widetilde{{\cal{R}}}_t^g({\cal{Z}}_{n\times d}):=\{\eta_t:\eta_\cdot\in\widetilde{{\cal{R}}}^g({\cal{Z}}_{n\times d})\}$ and when $n=1$, $L^p({\mathcal{F}}_t):=L_1^p({\mathcal{F}}_t)$, $H^p:=H_{1\times d}^p$, and employ a similar treatment for other spaces.\\

We observe that the BSDE$(g,T,\xi)$ always has a solution for some terminal variables. In particular, we have the following characterization for the existence of BSDEs.\\\\
\textbf{Proposition 2.3} \emph{For $g\in{\cal{G}}_n$ and $\xi\in L_n^0({\mathcal{F}}_T)$, the BSDE$(g,T,\xi)$ has a solution if and only if $\xi\in\widetilde{{\cal{R}}}_T^g(H_{n\times d}^2)$}.\\\\
\emph{Proof.} If $\xi\in\widetilde{{\cal{R}}}_T^g(H_{n\times d}^2)$, then there exist $y\in{\mathbf{R}}^n$ and $Z_t\in H_{n\times d}^2$ such that
$$Y^{y,Z,g}_t=y-\int_{0}^tg(s,Y^{y,Z,g}_s,Z_s)ds+\int_0^tZ_sdB_s\ \ \textmd{with} \ Y^{y,Z,g}_T=\xi.$$
It follows that $(Y^{y,Z,g}_t,Z_t)$ is a solution to the BSDE$(g,T,\xi)$. Conversely, if the BSDE$(g,T,\xi)$ has a solution $(Y_t,Z_t)\in{\mathcal{S}}_n\times H_{n\times d}^2$, then by Lemma 2.2, we have $\xi\in\widetilde{{\cal{R}}}_T^g(H_{n\times d}^2)$. \ \ $\Box$\\\\
\textbf{Remark 2.4} \emph{In Lemma 2.2, the initial data $(0,y)$ of the SDE(2.2) can be replaced by any $(\tau,\xi)\in({\cal{T}}_{0,T},L_n^0({\cal{F}}_\tau))$ and similarly, in Definition 2.1 and Proposition 2.3, the terminal data $(T,\xi)$ of the BSDE(2.1) can be replaced by any $(\tau,\xi)\in({\cal{T}}_{0,T},L_n^0({\cal{F}}_\tau))$.}\\

From Proposition 2.3, it follows that solving the BSDE$(g,T,\xi)$ is exactly to determine its solvable set $\widetilde{{\cal{R}}}_T^g(H_{n\times d}^2)$. Although it is difficult to recognise this set directly, it is possible that $\widetilde{{\cal{R}}}_T^g(H_{n\times d}^2)$ is computed by numerical methods. Now, using some known results, we compare $\widetilde{{\cal{R}}}_T^g(H_{n\times d}^2)$ with $L^p$ spaces in some cases.\\\\
\textbf{Proposition 2.5} \emph{(i) $\widetilde{{\cal{R}}}_T^g({{\cal{H}}_{n\times d}^p})=L_n^p({\mathcal {F}}_T)$ and $\widetilde{{\cal{R}}}^g({{\cal{H}}_{n\times d}^p})\subset {\cal{S}}_n^p$, when $p>1$, $g\in{\cal{L}}_n^\nu$ and $g(t,0,0)\in{\cal{H}}^p_{1\times1}$;}

\emph{(ii) $L^\infty({\mathcal {F}}_T)\subset\widetilde{{\cal{R}}}_T^g({{\cal{H}}^{BMO}})$, when $g\in{\cal{L}}^\rho\cap\widehat{{\cal{G}}}$ with $\rho(|\cdot|)=\nu(|\cdot|+1)$ and $|g(\cdot,0,0)|\leq k$, $k>0$;}

\emph{(iii) $\widetilde{{\cal{R}}}_T^g(H^2)\subset\widetilde{{\cal{R}}}_T^f(H^2)$, and $\xi\in\widetilde{{\cal{R}}}_T^f(H^2)$ if and only if $\exp(2\nu(\xi+\int_0^T\gamma_tdt))\in L^1({\mathcal {F}}_T)$, when $g,f\in{\cal{G}}^z$ such that $g(t,\cdot)\geq f(t,\cdot)=\nu|\cdot|^2+\gamma_t$, $dt\times dP$-a.e., where $\gamma_t\in H^1_{1\times1}$;}

\emph{(iv) ${\cal{R}}\subset\overline{{\cal{R}}}\subset\widetilde{{\cal{R}}}_T^g({\cal{H}}^\infty)\subset L^p({\mathcal {F}}_T),$ when $p>1$, $g\in{\cal{L}}^\rho\cap\widehat{{\cal{G}}}$ and $g(t,0,0)\in l^2_{1\times1}(0,T)$.}\\\\
\emph{Proof.} Proof of (i): Let $p>1$. It follows from Briand et al. \cite[Theorem 4.2]{BD} that $L_n^p({\mathcal {F}}_T)\subset\widetilde{{\cal{R}}}_T^g({{\cal{H}}_{n\times d}^p})$. Conversely, for $Y_t\in\widetilde{{\cal{R}}}^g({{\cal{H}}_{n\times d}^p})$, there exists a pair $(y,Z_t)\in\textbf{R}^n\times{\cal{H}}_{n\times d}^p$ such that $$Y_t=y-\int_{0}^tg(s,Y_s,Z_s)ds+\int_0^tZ_sdB_s.$$ Set $\tau_m:=\inf\{t\geq0:|Y_t|\geq m\}\wedge T$ and $Y_t^m:=Y_{t\wedge{\tau_m}},\ m\geq1.$ By BDG inequality and H\"older's inequality, we have
\begin{eqnarray*}
\ \ \ \ \ \ \ E\left[\sup_{r\in[0,T]}|Y_r^m|^p\right]&\leq&3^{p-1}|y|^p+3^{p-1}E\left(\int_{0}^T(\mu|Y_s^m|+\nu|Z_s|+|g(s,0,0)|)ds\right)^p\\
&&\ \ \ +3^{p-1}E\left(\sup_{r\in[0,T]}\left|\int_0^r Z_sdB_s\right|^p\right)\\
&\leq&C\left(1+E\left(\int_0^T\sup_{r\in[0,s]}|Y_r^m|^pds\right)+E\left[\left(\int_0^T|Z_s|^2ds\right)^\frac{p}{2}\right]\right),\ \ \ \ \ \ \ (2.3)
\end{eqnarray*}
where $C>0$ is a constant depending only on $y$, $\mu$, $\nu$, $T$, $p$ and $\|g(t,0,0)\|_{{\cal{H}}^p}$. Then by Gronwall inequality and letting $m\rightarrow\infty,$ we get $Y_\cdot\in{\cal{S}}_n^p$. Thus, $\widetilde{{\cal{R}}}^g({{\cal{H}}_{n\times d}^p})\subset {\cal{S}}_n^p$ and $\widetilde{{\cal{R}}}_T^g({{\cal{H}}_{n\times d}^p})\subset L_n^p({\mathcal {F}}_T)$.

Proof of (ii): (ii) can be derived from \cite[Theorem 2.2]{BE}.

Proof of (iii): For $\xi\in\widetilde{{\cal{R}}}_T^f(H^2)$, there exists a pair $(y,Z_t)\in\textbf{R}\times H^2$ such that $$\xi+\int_0^T\gamma_sds=y-\int_{0}^T\nu|Z_s|^2ds+\int_0^T Z_sdB_s.$$
By Zheng et al. \cite[Proposition 3.1]{Zheng3}, we can get $\exp(2\nu(\xi+\int_0^T\gamma_sds))\in L^1({\mathcal {F}}_T)$. Conversely, for $\xi$ such that $\exp(2\nu(\xi+\int_0^T\gamma_sds))\in L^1({\mathcal {F}}_T),$ by martingale representation theorem (see \cite[Theorem 2.46]{PR}) and setting $y_t:=E[\exp(2\nu(\xi+\int_0^T\gamma_sds))|{\cal{F}}_t]$, there exists $z_s\in H^2$ such that
$$y_t=\exp(2\nu(\xi+\int_0^T\gamma_s ds))-\int_t^Tz_sdB_s,\ \ t\in[0,T].$$
Moreover, by Doob's optional stopping theorem, we get that for each $\tau\in{\cal{T}}_{0,T}$, we have
  $$y_\tau=E[\exp(2\nu(\xi+\int_0^T\gamma_sds))|{\cal{F}}_\tau],$$
which implies that for each $\tau\in{\cal{T}}_{0,T}$, $y_\tau>0$. Set $$\tau:=\inf\{t\geq0:y_t\leq0\}\wedge T.$$ Since for each $\tau\in{\cal{T}}_{0,T}$, $y_\tau>0$, by the continuity of $y_t$, we can deduce $P(\tau=T)=1,$ which implies
$$P(y_t>0,t\in[0,T])=P(\tau=T)=1.$$
By applying It\^{o}'s formula to $\frac{1}{2\nu}\ln(y_t)$, we get that there exists $Z_t\in H^2$ such that
$$\xi=\frac{1}{2\nu}\ln(y_0)-\int_{0}^T(\nu|Z_s|^2+\gamma_s)ds+\int_0^T Z_sdB_s\ \  \textmd{with}\ \ Z_t=\frac{z_t}{2\nu y_t},$$
which implies that $\xi\in\widetilde{{\cal{R}}}_T^f(H^2).$ Thus, $\xi\in\widetilde{{\cal{R}}}_T^f(H^2)$ if and only if $\exp(2\nu(\xi+\int_0^T\gamma_sds))\in L^1({\mathcal {F}}_T).$ Moreover, since $g\geq f$, for $\xi\in\widetilde{{\cal{R}}}_T^g(H^2)$, there exists a pair $(y,Z_t)\in\textbf{R}\times H^2$ such that
$$\xi=y-\int_{0}^Tg(s,Z_s)ds+\int_0^T Z_sdB_s\leq y-\int_{0}^Tf(s,Z_s)ds+\int_0^T Z_sdB_s,$$
which implies that $\exp(2\nu(\xi+\int_0^T\gamma_sds))\in L^1({\mathcal {F}}_T).$ Thus, (iii) holds true.

Proof of (iv): Clearly, ${\cal{R}}\subset\overline{{\cal{R}}}$. Similar to (2.3), we can get $\widetilde{{\cal{R}}}_T^g({\cal{H}}^\infty)\subset L^p({\mathcal {F}}_T),\ p>1.$ Moreover, by \cite[Theorem 2.2]{CN}, we have $\overline{{\cal{R}}}\subset\widetilde{{\cal{R}}}_T^g({\cal{H}}^\infty)$. Thus, (iv) holds true.\ \ $\Box$\\\\
\textbf{Remark 2.6} \emph{In Proposition 2.5(ii), in general, $L^\infty({\mathcal {F}}_T)\neq\widetilde{{\cal{R}}}_T^g({{\cal{H}}^{BMO}}),$ since for $(y,z)\in{\bf{R}}^n\times{\bf{R}}^{n\times d}$ and $g\equiv0$, $y+zB_T$ is unbounded. In Proposition 2.5(iv), in general, $\overline{{\cal{R}}}\neq\widetilde{{\cal{R}}}_T^g({\cal{H}}^\infty)$, since for $Z_s\in{\cal{H}}^{\infty}$ which is not Malliavin differentiable, $y-\int_0^Tg(t,Z_t)dt+\int_0^TZ_tdB_t$ may be not Malliavin differentiable (see \cite[Lemma 1.3.4]{Na}). According to Proposition 2.5(iii), when $g\in{\cal{G}}^z$ has a superquadratic growth in $z$ such that $g(t,\cdot)\geq\nu|\cdot|^2+\gamma_t,$ $dt\times dP-a.e.,$ with $\gamma_t\in H^1_{1\times1}$, to guarantee the BSDE$(g,T,\xi)$ has a solution, $\xi$ at least satisfies the integrability: $\exp(2\nu(\xi+\int_0^T\gamma_tdt))\in L^1({\mathcal {F}}_T)$.}\\

It is well-known that BSDEs may have multiple solutions, even for the case that the terminal variable is $L^2$ integrable and $g\equiv0$ (see Xing \cite{Xing}). This implies that when we discuss the uniqueness of solution, the space of solution must be emphasized. For $g\in{\cal{G}}_n$ and $\xi\in\widetilde{{\cal{R}}}_T^g(H_{n\times d}^2)$, we can deduce that for two solutions $(Y_t,Z_t)$ and $(y_t,z_t)$ of the BSDE$(g,T,\xi)$,
$$(Y_\cdot,Z_\cdot)=(y_\cdot,z_\cdot)\Longleftrightarrow Y_\cdot=y_\cdot\Longleftrightarrow Z_\cdot=z_\cdot,\ \ dt\times dP-a.e.,$$
where the first equivalent is clear, and the second equivalent can be derived from
$$|Y_t-y_t|\leq\int_t^T\mu|Y_s-y_s|ds,\ \ \  (\textmd{by the Lipschitz condition on}\ y)$$
and backward Gronwall inequality (\cite[Corollary 6.62]{PR}).

Let ${\cal{Y}}_n$ denote a subset of ${\cal{S}}_n$ satisfying (i) and (ii):

(i) for each $\eta_t,\zeta_t\in{\cal{Y}}_n$ and $\tau\in{\cal{T}}_{0,T}$ with $\eta_\tau=\zeta_\tau$, we have $\eta_t1_{[0,\tau)}(t)+\zeta_t1_{[\tau,T]}(t)\in{\cal{Y}}_n;$

(ii) for $y_t\in{\cal{S}}_n$, we have $y_t\in{\cal{Y}}_n,$ if there exist $\eta_t^1,\eta_t^2,\zeta_t^1,\zeta_t^2\in{\cal{Y}}_n$ and $\tau\in{\cal{T}}_{0,T}$, such that $$\eta_t^11_{[0,\tau)}(t)+\eta_t^21_{[\tau,T]}(t)\leq y_t\leq\zeta_t^11_{[0,\tau)}(t)+\zeta_t^21_{[\tau,T]}(t),\ \ \forall t\in[0,T].$$
Clearly, ${\cal{S}}_n$ and ${\cal{S}}_n^p,\ p\geq1,$ are some examples of ${\cal{Y}}_n$. For convenience, for a space of stochastic process $\cal{O}$ and $\tau,\sigma\in{\cal{T}}_{0,T}$ such that $\tau\leq\sigma,$ we set ${\cal{O}}^{[\tau,\sigma]}:=\{(x_t)_{t\in[\tau,\sigma]}: x_t\in{\cal{O}}\}$.

For $g\in{\cal{G}}_n$ and $\xi\in L_n^0({\cal{F}}_T)$, we say the BSDE$(g,T,\xi)$ has a \textbf{unique} (resp. \textbf{maximal}) \textbf{solution} $(Y_t,Z_t)\in{\cal{Y}}_n\times{\cal{Z}}_{n\times d}$, if for any $\tau\in{\cal{T}}_{0,T}$ and any solution $(y_t,z_t)\in{\cal{Y}}_n^{[\tau,T]}\times{\cal{Z}}_{n\times d}^{[\tau,T]}$ to the BSDE$(g,T,\xi)$ on the interval $[\tau,T]$, we have $Y_t=y_t$ (resp. $Y_t\geq y_t$) for each $t\in[\tau,T].$\\

When $n=1$, we have the following general comparison theorems. \\\\
\textbf{Theorem 2.7} \emph{(i) Let} $g\in\textbf{G}$ \emph{and $\xi, \eta\in L^0({\cal{F}}_T)$ such that $\xi\geq\eta$. If the BSDE$(g,T,\xi)$ has a unique (resp. maximal) solution $(Y_t,Z_t)\in{\cal{Y}}\times{\cal{Z}}$ and the BSDE$(g,T,\eta)$ has a solution $(y_t,z_t)\in{\cal{Y}}\times{\cal{Z}}$, then for each $t\in[0,T]$, $Y_t\geq y_t$.}

\emph{(ii) Let $g\in{\cal{G}}$, $f_t\in H^1_{1\times1}$ and $\xi, \eta\in L^0({\cal{F}}_T)$ such that $\xi\geq\eta$. If the BSDE$(g,T,\xi)$ have a unique (resp. maximal) solution $(Y_t,Z_t)\in{\cal{Y}}\times{\cal{Z}}$ and the BSDE$(f_t,T,\xi)$ have a solution $(y_t,z_t)\in{\cal{Y}}\times{\cal{Z}}$ such that $g(t,y_t,z_t)\geq f_t$, $dt\times dP$-$a.e.,$ then for each $t\in[0,T]$, $Y_t\geq y_t$.}
\\\\
\emph{Proof.} In the case (i) (or (ii)), if there exists $t_0\in[0,T)$ such that
$$P(\{Y_{t_0}<y_{t_0}\})>0,$$
then by the fact
$$\bigcup_{m\geq1}\{Y_{t_0}\leq y_{t_0}-\frac{1}{m}\}=\{Y_{t_0}<y_{t_0}\},$$
we get that there exists $l\geq1$ such that
$$P(\{Y_{t_0}\leq y_{t_0}-\frac{1}{l}\})>0.\eqno(2.3)$$
We define the following stopping times:
$$\sigma_l:=\inf\{t\geq t_0:Y_t\leq y_t-\frac{1}{l}\}\wedge T,\eqno(2.4)$$
and
$$\delta:=\inf\{t\geq\sigma_l:Y_t\geq y_t\}\wedge T.$$
By the continuity of $Y_t$ and $y_t$, we can conclude that
$$Y_{\sigma_l}\leq y_{\sigma_l}-\frac{1}{l}\ \ \ \textmd{on the set}\ \ \{\sigma_l<T\},\eqno(2.5)$$
and $$\{\sigma_l<\delta\}=\{\sigma_l<T\}.\eqno(2.6)$$
By (2.5), (2.6), the fact $\xi\geq\eta$, and the continuity of $Y_t$ and $y_t$, we can get
$$Y_\delta=y_\delta \ \ \ \textmd{on the set}\ \{\sigma_l<\delta\}.\eqno(2.7)$$
By (2.4) and (2.6), we have
$$\{Y_{t_0}\leq y_{t_0}-\frac{1}{l}\}\subset\{\sigma_l<\delta\},$$
which together with (2.3) gives
$$P(\{\sigma_l<\delta\})>0.\eqno(2.8)$$

\textbf{Proof of (i):}  Set
$$(\tilde{y}_t, \tilde{z}_t):=(y_t1_{[\sigma_l,\delta)}(t)+Y_t1_{[\delta,T]}(t), \ z_t1_{[\sigma_l,\delta)}(t)+Z_t1_{[\delta,T]}(t)),\ \ t\in[\sigma_l,T].\eqno(2.9)$$
Then we have
$$\tilde{y}_{\sigma_l}=y_{\sigma_l}\ \ \ \ \textmd{on the set}\ \ \{\sigma_l<\delta\},\eqno(2.10)$$
Since
$$y_{t\vee\sigma_l}=y_\delta+\int_{t\vee\sigma_l}^\delta g(s,y_s,z_s)ds+\int_{t\vee\sigma_l}^\delta z_sdB_s, \ t\in[0,\delta],$$
we have
$$y_{t\vee\sigma_l}1_{\{\sigma_l<\delta\}}+y_\delta1_{\{\sigma_l=\delta\}}=y_\delta1_{\{\sigma_l<\delta\}}+y_\delta1_{\{\sigma_l=\delta\}}+\int_{t\vee\sigma_l}^\delta g(s,y_s,z_s)ds+\int_{t\vee\sigma_l}^\delta z_sdB_s, \ t\in[0,\delta],$$
which can gives
\begin{eqnarray*}
&&y_{t\vee\sigma_l}1_{\{\sigma_l<\delta\}}+Y_\delta1_{\{\sigma_l=\delta\}}\\&=&y_\delta1_{\{\sigma_l<\delta\}}+Y_\delta1_{\{\sigma_l=\delta\}}+\int_{t\vee\sigma_l}^\delta g(s,y_s,z_s)ds+\int_{t\vee\sigma_l}^\delta z_sdB_s\\
&=&Y_\delta+\int_{t\vee\sigma_l}^\delta g(s,y_s,z_s)ds+\int_{t\vee\sigma_l}^\delta z_sdB_s \ \ \ \ \ (\textmd{by}\ (2.7))\\
&=&\tilde{y}_\delta+\int_{t\vee\sigma_l}^\delta g(s,\tilde{y}_s,\tilde{z}_s)ds+\int_{t\vee\sigma_l}^\delta \tilde{z}_sdB_s, \ \ t\in[0,\delta]. \ \ \ \ \ (\textmd{by}\ (2.9))
\end{eqnarray*}
This together with
$$\tilde{y}_{t\vee\sigma_l}=y_{t\vee\sigma_l}1_{\{\sigma_l<\delta\}}+Y_\delta1_{\{\sigma_l=\delta\}}, \ \ t\in[0,\delta],\ \ \ \ \ (\textmd{by}\ (2.7)\ \textmd{and}\ (2.9))$$
implies
$$\tilde{y}_{t\vee\sigma_l}=\tilde{y}_\delta+\int_{t\vee\sigma_l}^\delta g(s,\tilde{y}_s,\tilde{z}_s)ds+\int_{t\vee\sigma_l}^\delta \tilde{z}_sdB_s, \ t\in[0,\delta].$$
By (2.9), we also have
$$\tilde{y}_{t\vee\delta}=\xi+\int_{t\vee\delta}^T g(s,\tilde{y}_s,\tilde{z}_s)ds+\int_{t\vee\delta}^T \tilde{z}_sdB_s, \ t\in[0,T].$$
The above two equations imply that $(\tilde{y}_t,\tilde{z}_t)\in{\cal{Y}}^{[\sigma_l,T]}\times{\cal{Z}}^{[\sigma_l,T]}$ is a solution to the BSDE$(g,T,\xi)$ on the interval $[\sigma_l,T]$. Moreover, since $(Y_t,Z_t)$ is a unique (resp. maximal) solution of the BSDE$(g,T,\xi)$, we get $Y_{\sigma_l}=\tilde{y}_{\sigma_l}$ (resp. $Y_{\sigma_l}\geq\tilde{y}_{\sigma_l}$), which together with (2.10) implies
$$Y_{\sigma_l}=y_{\sigma_l}\ (\textmd{resp.}\ Y_{\sigma_l}\geq y_{\sigma_l})\ \ \ \ \textmd{on the set}\ \ \{\sigma_l<\delta\}.$$
This together with (2.5) and (2.6) implies $P(\{\sigma_l<\delta\})=0$, which contradicts (2.8). Thus, (i) holds true.

\textbf{Proof of (ii):}  We consider the following two SDEs on the interval $[\sigma_l,T]$:
$$\tilde{Y}_{t\vee\sigma_l}=y_{\sigma_l}-\int_{\sigma_l}^{t\vee\sigma_l}g(s,\tilde{Y}_s,z_s)ds+\int_{\sigma_l}^{t\vee\sigma_l}z_sdB_s, \ t\in[0,T]$$ and $$y_{t\vee\sigma_l}=y_{\sigma_l}-\int_{\sigma_l}^{t\vee\sigma_l}f_sds+\int_{\sigma_l}^{t\vee\sigma_l}z_sdB_s, \ t\in[0,T].$$
Since $g(t,y_t,z_t)\geq f_t$, $dt\times dP$-$a.e.$, by the comparison result in Lemma 2.2, we get that for each $t\in[0,T]$, $$\tilde{Y}_{t\vee\sigma_l}\leq y_{t\vee\sigma_l}.\eqno(2.11)$$
We define the stopping time:
$$\kappa:=\inf\{t\geq\sigma_l:Y_t\geq\tilde{Y}_t\}\wedge T.$$
Then similar to (2.6), we have
$$\{\sigma_l<\kappa\}=\{\sigma_l<T\}.\eqno(2.12)$$
Since
$$\tilde{Y}_T\leq y_T=\eta\leq\xi=Y_T,\ \ \ (\textmd{by}\ (2.11))$$ by a similar argument as in (2.7), we get
$$\tilde{Y}_\kappa=Y_\kappa \ \ \textmd{on the set} \ \{\sigma_l<\kappa\}.
\eqno(2.13)$$
For $t\in[\sigma_l,T]$, we define
$$\breve{Y}_t:=\tilde{Y}_t1_{[\sigma_l,\kappa)}(t)+Y_t1_{[\kappa,T]}(t)$$
and
$$\hat{Y}_t:=y_t1_{[\sigma_l,\delta)}(t)+Y_t1_{[\delta,T]}(t).$$
Then by (2.7), (2.12) and (2.13), we get that $\breve{Y}_t\in{\cal{S}}^{[\sigma_l,T]}$ and $\hat{Y}_t\in{\cal{Y}}^{[\sigma_l,T]}$. Moreover, by (2.11) and the definitions of $\delta$ and $\kappa$, we can get $\kappa\leq \delta$, and for each $t\in[0,T]$, $$\breve{Y}_t1_{[\sigma_l,\kappa)}(t)=\tilde{Y}_t1_{[\sigma_l,\kappa)}(t)\geq Y_t1_{[\sigma_l,\kappa)}(t)\ \ \ \textmd{and}\ \ \ \breve{Y}_t1_{[\sigma_l,\delta)}(t)\leq y_t1_{[\sigma_l,\delta)}(t)=\hat{Y}_t1_{[\sigma_l,\delta)}(t).$$ This together with the fact that $\kappa\leq \delta$, implies that for each $t\in[\sigma_l,T]$, $${Y}_t\leq\breve{Y}_t\leq\hat{Y}_t.$$
Thus for each $t\in[0,T]$, $$y_t1_{[0,\sigma_l)}(t)+{Y}_t1_{[\sigma_l,T]}(t)\leq y_t1_{[0,\sigma_l)}(t)+\breve{Y}_t1_{[\sigma_l,T]}(t)\leq y_t1_{[0,\sigma_l)}(t)+\hat{Y}_t1_{[\sigma_l,T]}(t).$$ In view of $y_t1_{[0,\sigma_l)}(t)+\breve{Y}_t1_{[\sigma_l,T]}(t)\in{\cal{S}}$, by the definition of ${\cal{Y}}$, we have $$y_t1_{[0,\sigma_l)}(t)+\breve{Y}_t1_{[\sigma_l,T]}(t)\in{\cal{Y}}, \ \ \textmd{with}\ \ \breve{Y}_t1_{[\sigma_l,\kappa)}(t)=\tilde{Y}_t1_{[\sigma_l,\kappa)}(t),\ \forall t\in[0,T].\eqno(2.14)$$.

Set
$$(\bar{y}_t,\bar{z}_t):=(\tilde{Y}_t1_{[\sigma_l,\kappa)}(t)+Y_t1_{[\kappa,T]}(t),\ z_s1_{[\sigma_l,\kappa)}(t)+Z_t1_{[\kappa,T]}(t)),\ \ t\in[0,T].$$
Then by (2.6), (2.12) and the fact that $\tilde{Y}_{\sigma_l}=y_{\sigma_l}$, we have
$$\bar{y}_{\sigma_l}=y_{\sigma_l}\ \ \ \textmd{on the set} \ \ \{\sigma_l<\delta\},\eqno(2.15)$$
and by (2.13), the similar arguments as in (i), and (2.14), we can check that $(\bar{y}_t,\bar{z}_t)_{t\in[\sigma_l,T]}\in{\cal{Y}}^{[\sigma_l,T]}\times{\cal{Z}}^{[\sigma_l,T]}$ is a solution to the BSDE$(g,T,\xi)$ on the interval $[\sigma_l,T]$. Since $(Y_t,Z_t)$ is a unique (resp. maximal) solution of the BSDE$(g,T,\xi)$, we have $Y_{\sigma_l}=\bar{y}_{\sigma_l}$ (resp. $Y_{\sigma_l}\geq\bar{y}_{\sigma_l}$), which together with (2.15) implies
$$Y_{\sigma_l}=y_{\sigma_l}\ (\textmd{resp.}\ Y_{\sigma_l}\geq y_{\sigma_l})\ \ \ \textmd{on the set} \ \ \{\sigma_l<\delta\}.$$
This together with (2.5) and (2.6) implies $P(\{\sigma_l<\delta\})=0$, which contradicts (2.8). Thus, (ii) holds true.

The proof is complete.\ \ $\Box$\\

The proofs of Theorem 2.7(i)(ii) rely mainly on the uniqueness of solution. Theorem 2.7(i) is proved without using any assumption on $(y,z)$. Theorem 2.7(ii) is proved only by using the Lipschitz assumption on $y$ (without using any assumption on $z$), in order to use the comparison theorem for a special SDE. In fact, this Lipschitz assumption can be further generalized, since the comparison theorem for SDEs holds true under some general conditions such like the monotonic condition (see \cite[Proposition 3.14]{PR}), etc. Consequently, Theorem 2.7 holds true for the BSDEs appearing in \cite{EPQ}, \cite{K}, \cite{BD}, \cite{BH}, \cite{Jia}, \cite{BE} and \cite{CN}, etc. This means that comparison theorem usually holds true as long as the BSDE has a unique (resp. maximal) solution in some space.

The following Proposition 2.8 gives an explicit solution to a special $n$-dimensional mean-field BSDE with terminal variables set $\widetilde{{\cal{R}}}_T^0(l_{n\times d}^2(0,T))$, from which an invariant representation theorem is derived without any assumption on $z$.
\\\\
\textbf{Proposition 2.8} \emph{Let $(t,y,z)\in[0,T)\times{\mathbf{R}}^n\times{\mathbf{R}}^{n\times d}$ and $\varepsilon\in[0,T-t]$. Then for $g\in\widehat{{\cal{G}}}_n$ and $h(s)\in l_{n\times d}^2(0,t+\varepsilon)$ such that $\int_0^{t+\varepsilon}|g(s,0,h(s))|ds<\infty$, the following mean-field BSDE
$$Y^{t+\varepsilon}_s=y+\int_0^{t+\varepsilon} h(r)dB_r+\int_{s}^{t+\varepsilon} g(r,E[Y^{t+\varepsilon}_r],Z^{t+\varepsilon}_r)
dr-\int_s^{t+\varepsilon}Z^{t+\varepsilon}_rdB_r,\ \ s\in[0,t+\varepsilon], \eqno(2.16)$$
has an adapt solution $$(Y^{t+\varepsilon}_s,Z^{t+\varepsilon}_s)=\left(\varphi^{t+\varepsilon}(s)+\int_0^sh(r)dB_r,\ h(s)\right),\ s\in[0,{t+\varepsilon}],\eqno(2.17)$$
where $\varphi^{t+\varepsilon}(s)\in l_{n\times d}^2(0,t+\varepsilon)$ is a solution to the backward ODE: $$\varphi^{t+\varepsilon}(s)=y+\int_s^{t+\varepsilon}g(r,\varphi^{t+\varepsilon}(r),h(r))dr,\ \ s\in[0,t+\varepsilon].$$ Moreover, if $h(s)=z 1_{s\in[t,t+\varepsilon]}(s)$, $s\in[0,t+\varepsilon]$, then for almost each $t\in[0,T),$ we have the following invariant representation:} $$g(t,y,z)=\lim_{\varepsilon\rightarrow0}\frac{1}{\varepsilon}(Y^{t+\varepsilon}_t-y).\eqno(2.18)$$
\emph{Proof.} It follows from the ODE theory that
$$\varphi^{t+\varepsilon}(s)=y+\int_s^{t+\varepsilon}g(r,\varphi^{t+\varepsilon}(r),h(r))dr,\ s\in[0,t+\varepsilon],$$
has a unique solution $\varphi^{t+\varepsilon}(s)\in l_{n\times d}^2(0,t+\varepsilon)$ (This result can also be got from \cite[Theorem 5.17]{PR}). In view of $E[Y^{t+\varepsilon}_s]=\varphi^{t+\varepsilon}(s)$, it is not hard to check that (2.17) is a solution to the BSDE(2.16). Moreover, if $h(s)=z 1_{s\in[t,t+\varepsilon]}(s)$, $s\in[0,t+\varepsilon]$, then for each $s\in[t,t+\varepsilon],$ we have
\begin{eqnarray*}
|\varphi^{t+\varepsilon}(s)-y|
&=&\left|\int_s^{t+\varepsilon}\textmd(g(r,(\varphi^{t+\varepsilon}(r)-y)+y,z)-g(r,y,z))+g(r,y,z))dr\right|\\
&\leq&\int_s^{t+\varepsilon}\mu|\varphi^{t+\varepsilon}(r)-y|dr+\int_s^{t+\varepsilon}|g(r,y,z)|dr
\end{eqnarray*}
By backward Gronwall inequality, we have for each $s\in[t,t+\varepsilon],$
$$|\varphi^{t+\varepsilon}(s)-y|\leq e^{\mu(t+\varepsilon-s)}\int_s^{t+\varepsilon}|g(r,y,z)|dr\leq e^{\mu\varepsilon}\int_t^{t+\varepsilon}|g(r,y,z)|dr.\eqno(2.19)$$
By (2.17) and (2.19), we have
\begin{eqnarray*}
\left|\frac{1}{\varepsilon}(Y^{t+\varepsilon}_t-y)-g(t,y,z)\right|&=&\left|\frac{1}{\varepsilon}\int_t^{t+\varepsilon}g(r,\varphi^{t+\varepsilon}(r),z)dr-g(t,y,z)\right|\\
&\leq&\frac{1}{\varepsilon}\int_t^{t+\varepsilon}(|g(r,\varphi^{t+\varepsilon}(r),z)-g(r,y,z)|+|g(r,y,z)-g(t,y,z)|)dr\\
&\leq&\frac{1}{\varepsilon}\int_t^{t+\varepsilon}\mu|\varphi^{t+\varepsilon}(r)-y|dr+\frac{1}{\varepsilon}\int_t^{t+\varepsilon}|g(r,y,z)-g(t,y,z)|dr\\
&\leq&\mu e^{\mu\varepsilon}\int_t^{t+\varepsilon}|g(u,y,z)|du+\frac{1}{\varepsilon}\int_t^{t+\varepsilon}|g(r,y,z)-g(t,y,z)|dr.
\end{eqnarray*}
Then by the continuity of integration and Lebesgue's Lemma, we obtain (2.18). \ \ $\Box$\\\\
\textbf{Example 2.9}  Let $\xi\in L_n^2({\cal{F}}_T)$. By martingale representation theorem, there exists a unique $Z_s\in{\cal{H}}_{n\times d}^2$ such that $\xi=E[\xi]+\int_0^TZ_sdB_s$. Then similar to Proposition 2.8, we have

(i) For $g\in\widehat{{\cal{G}}}_n$, when $\int_0^T|g(s,0,E[Z_s])|ds<\infty,$ the mean-field BSDE
$$Y_t=\xi+\int_{t}^Tg(s,E[Y_s],E[Z_s])
ds-\int_t^T Z_sdB_s,\ \ \ t\in[0,T],\eqno(2.20)$$
has an adapt solution
$(Y_t,Z_t)=(\varphi(t)+\int_0^tZ_sdB_s, Z_t), t\in[0,T],$ where $\varphi(t)\in l_{n\times d}^2(0,T)$ is a solution to the backward ODE: $\varphi(t)=E[\xi]+\int_t^Tg(s,\varphi(s),E[Z_s])ds,\ t\in[0,T].$

(ii) For $g\in{\cal{G}}_n$, when $\int_{0}^TE|g(s,0,Z_s)|ds<\infty$, the mean-field BSDE
$$Y_t=\xi+\int_{t}^TE[g(s,E[Y_s],Z_s)]ds-\int_t^T Z_sdB_s,\ \ t\in[0,T],$$
has an adapt solution
$(Y_t,Z_t)=(\varphi(t)+\int_0^tZ_sdB_s, Z_t), t\in[0,T],$ where $\varphi(t)\in l_{n\times d}^2(0,T)$ is a solution to the backward ODE: $\varphi(t)=E[\xi]+\int_{t}^TE[g(s,\varphi(s),Z_s)]ds,\ t\in[0,T].$\\\\
\textbf{Remark 2.10} \emph{Mean-field BSDEs with general structures have been studied by Buckdahn et al. \cite{Li} for $L_n^2$ terminal variables under Lipschitz condition, and by Hibon et al. \cite{HH} for $L^\infty$ terminal variables under quadratic growth condition. In general, the comparison theorem is not true for mean-field BSDEs (see \cite{Li}). But when $n=1$, $g$ is independent of $z$, for $\xi,\eta\in L^2$ such that $E[\xi]\geq E[\eta]$, we can get that for each $t\in[0,T]$, $E[Y_t^1]\geq E[Y_t^2]$, where $(Y_t^1,Z_t^1)$ and $(Y_t^2,Z_t^2)$ are the solutions to the BSDE(2.20) with terminal variables $\xi$ and $\eta$, respectively, given as in Example 2.9(i).}\\

From the invariant representation in Proposition 2.8, we can get the following $n$-dimensional converse comparison result by passing a limit of a series of solutions directly. We omit its proof.\\\\
\textbf{Proposition 2.11} \emph{For $g, f\in\widehat{{\cal{G}}}_n\cap{{\cal{G}}}_n^z$, $(t,z)\in[0,T)\times{\mathbf{R}}^{n\times d}$ and $\varepsilon\in[0,T-t]$, let $(y_s^{t+\varepsilon,z}, z_s^{t+\varepsilon,z})$ and $(\tilde{y}_s^{t+\varepsilon,z}, \tilde{z}_s^{t+\varepsilon,z})$ be the solutions to the BSDE$(g,t+\varepsilon,z(B_{t+\varepsilon}-B_{t}))$ and the BSDE$(f,t+\varepsilon,z(B_{t+\varepsilon}-B_{t}))$ given as in (2.17), respectively. If for each $(t,z)\in[0,T)\times{\mathbf{R}}^{n\times d}$ and $\varepsilon\in[0,T-t]$, we have $y_t^{t+\varepsilon,z}\geq \tilde{y}_t^{t+\varepsilon,z},$ then for each $z\in{\mathbf{R}}^{n\times d},$ $g(t,z)\geq f(t,z),$ $a.e.,$ $t\in[0,T)$.}\\

At the end of this section, we give a uniqueness result, which generalizes the sets of terminal variables in \cite[Theorem 2.2 and Corollary 2.2]{BE} and \cite[Theorem 2.2]{CN}.\\\\
\textbf{Proposition 2.12} \emph{For $g\in{\cal{L}}^{\rho}$ with $g(t,0,0)\in{\cal{H}}^{BMO}_{1\times1}$, and $\xi\in\widetilde{{\cal{R}}}_T^g({\cal{H}}_\rho^{BMO})$, the BSDE$(g,T,\xi)$ has a unique solution in ${\cal{S}}\times{\cal{H}}_\rho^{BMO}$, in particular in $\bigcap_{p>1}{\cal{S}}^p\times{\cal{H}}_\rho^{BMO}$.}\\\\
\emph{Proof.} By an analogy to Proposition 2.3, we only need to prove the uniqueness. Let $(Y_t,Z_t)\in{\cal{S}}\times{\cal{H}}_\rho^{BMO}$ and $(\tilde{Y}_t,\tilde{Z}_t)\in{\cal{S}}\times{\cal{H}}_\rho^{BMO}$ be the solutions to the BSDE$(g,T,\xi)$. Since $Z_t, \tilde{Z}_t\in{\cal{H}}_\rho^{BMO}$ and
$$\rho(|Z_t|\vee|\tilde{Z}_t|)=\rho(|Z_t|)\vee\rho(|\tilde{Z}_t|)\leq\rho(|Z_t|)+\rho(|\tilde{Z}_t|),\ \ \ \forall t\in[0,T],$$
we get $$\sup_{\tau\in{\cal{T}}_{0,T}}\left\|E[\int_\tau^T\rho(|Z_t|\vee|\tilde{Z}_t|)^2dt|{\cal{F}}_\tau]\right\|_\infty<\infty.\eqno(2.21)$$

By a classical linearization argument, we have
$$Y_t-\tilde{Y}_t=\int_t^T(a_s(Y_s-\tilde{Y}_s)+b_s(Z_s-\tilde{Z}_s))ds-\int_t^T(Z_s-\tilde{Z}_s)dB_s,\ \ t\in[0,T],\eqno(2.22)$$
where
$$a_s:=\frac{g(s,Y_s,Z_s)-g(s,\tilde{Y}_s,Z_s)}{Y_s-\tilde{Y}_s}1_{|Y_s-\tilde{Y}_s|>0}$$
and $$b_s:=\frac{g(s,\tilde{Y}_s,Z_s)-g(s,\tilde{Y}_s,\tilde{Z}_s)}{|Z_s-\tilde{Z}_s|^2}(Z_s-\tilde{Z}_s)1_{|Z_s-\tilde{Z}_s|>0}.$$
Clearly, for each $s\in[0,T]$, we have $|a_s|\leq\mu$ and $|b_s|\leq\rho(|Z_s|\vee|\tilde{Z}_s|).$ By (2.21), we have $b_s\in{\cal{H}}^{BMO}.$ Thus by Girsanov theorem, there exists an equivalent probability $Q$, under which $B^b_t:=B_t-\int_0^tb_sds$ is standard Brownian motion. By (2.22) and It\^{o}'s formula, we can get
$$e^{\int_0^ta_rdr}(Y_t-\tilde{Y}_t)=-\int_t^Te^{\int_0^sa_rdr}(Z_s-\tilde{Z}_s)dB^b_s,\ \ t\in[0,T].\eqno(2.23)$$
Since $Z_t, \tilde{Z}_t\in{\cal{H}}^{BMO}$ under the probability $P$, we get that $Z_t, \tilde{Z}_t\in{\cal{H}}^{BMO}$ under the probability $Q$ (see \cite[Theorem 3.6]{Ka}). Hence, (2.23) implies that for each $t\in[0,T]$, $Y_t=\tilde{Y}_t$. This gives the uniqueness of solution.

Now, for the solution $(Y_t,Z_t)\in{\cal{S}}\times{\cal{H}}_\rho^{BMO}$ to the BSDE$(g,T,\xi)$, we will show $(Y_t,Z_t)\in\bigcap_{p>1}{\cal{S}}^p\times{\cal{H}}_\rho^{BMO}$. Since $Z_t\in{\cal{H}}_\rho^{BMO}$, by the energy inequality for BMO martingales (see \cite[Page 29]{Ka}), we have $Z_t\in\bigcap_{p>1}{\cal{H}}^p$ and for each $p>1$, $E[(\int_0^T\rho(|Z_t|)^2dt)^{\frac{p}{2}}]<\infty$. Then by a similar argument as in (2.3), we can get that $Y_t\in{\cal{S}}^p$ for each $p>1.$\ \ $\Box$
\section{$g$-expectations}
Using the results obtained in Section 2, we can further study the theory of $g$-expectations. Since the uniqueness of solutions to BSDEs depends on the space of solutions, we can not provide a unified definition of $g$-expectations for all the cases. We consider the following cases.\\\\
\textbf{Definition 3.1} (i) For $g\in{\cal{L}}_n^\nu\cap{\cal{G}}_n^0$ and $\xi\in\widetilde{{\cal{R}}}_T^g({\cal{H}}_{n\times d}^p)$, $p>1$, if $(Y_t,Z_t)$ is the solution to the BSDE$(g,T,\xi)$ such that $Z_t\in{\cal{H}}_{n\times d}^p$, then the conditional $g$-expectation of $\xi$ is defined by ${\cal{E}}_n^g[\xi|{\cal{F}}_t]:=Y_t$ for $t\in[0,T]$ and the $g$-expectation of $\xi$ is defined by ${\cal{E}}_n^g[\xi]:=Y_0$.

(ii) For $g\in{\cal{L}}^\rho\cap\widehat{{\cal{G}}}\cap{\cal{G}}^0$ with $\rho(|\cdot|)=\nu(|\cdot|+1)$ and $\xi\in\widetilde{{\cal{R}}}_T^g({\cal{H}}^{BMO})$, if $(Y_t,Z_t)$ is the solution to the BSDE$(g,T,\xi)$ such that $Z_t\in{\cal{H}}^{BMO}$, then the conditional $g$-expectation of $\xi$ is defined by ${\cal{E}}^g[\xi|{\cal{F}}_t]:=Y_t$ for $t\in[0,T]$ and the $g$-expectation of $\xi$ is defined by ${\cal{E}}^g[\xi]:=Y_0$.

(iii) For $g\in{\cal{L}}^\rho\cap\widehat{{\cal{G}}}\cap{\cal{G}}^0$ and $\xi\in\widetilde{{\cal{R}}}_T^g({\cal{H}}^\infty)$, if $(Y_t,Z_t)$ is the solution to the BSDE$(g,T,\xi)$ such that $Z_t\in{\cal{H}}^{\infty}$, then the conditional $g$-expectation of $\xi$ is defined by ${\cal{E}}^g[\xi|{\cal{F}}_t]:=Y_t$ for $t\in[0,T]$ and the $g$-expectation of $\xi$ is defined by ${\cal{E}}^g[\xi]:=Y_0$.

(iv) For $g\in\widehat{{\cal{G}}}_n\cap{\cal{G}}_n^z\cap{\cal{G}}_n^0$ and $\xi\in\widetilde{{\cal{R}}}_T^g(l^\infty_{n\times d}(0,T))$, if $(Y_t,Z_t)$ is the solution to the BSDE$(g,T,\xi)$ such that $Z_t\in l^\infty_{n\times d}(0,T)$, then the conditional $g$-expectation of $\xi$ is defined by ${\cal{E}}_n^g[\xi|{\cal{F}}_t]:=Y_t$ for $t\in[0,T]$ and the $g$-expectation of $\xi$ is defined by ${\cal{E}}_n^g[\xi]:=Y_0$. \\\\
\textbf{Remark 3.2} \emph{By \cite{BD}, Proposition 2.5 and Proposition 2.12, we conclude that Definition 3.1(i)-(iii) are well-defined, and in Definition 3.1(i), ${\cal{E}}_n^g[\xi|{\cal{F}}_t]\in{\cal{S}}_n^p$ and $\widetilde{{\cal{R}}}_T^g({\cal{H}}_{n\times d}^p)=L^p_n({\cal{F}}_T)$, in Definition 3.1(ii), ${\cal{E}}^g[\xi|{\cal{F}}_t]\in\bigcap_{p>1}{\cal{S}}^p$ and $L^\infty({\cal{F}}_T)\subset\widetilde{{\cal{R}}}_T^g({\cal{H}}^{BMO})$, and in Definition 3.1(iii), ${\cal{E}}^g[\xi|{\cal{F}}_t]\in\bigcap_{p>1}{\cal{S}}^p$ and $\overline{{\cal{R}}}\subset{{\cal{R}}}_T^g({\cal{H}}^\infty)$. We can also show that Definition 3.1(iv) is well-defined. In fact, for two solutions $(y_t,z_t)$ and $(\tilde{y}_t,\tilde{z}_t)$ to the BSDE($g,T,\xi$) with $z_t,\tilde{z}_t\in l^\infty_{n\times d}$, we have
$$y_0-\int_0^Tg(s,z_s)ds+\int_0^Tz_sdB_s=\tilde{y}_0-\int_0^Tg(s,\tilde{z}_s)ds+\int_0^T\tilde{z}_sdB_s.$$
By taking expectation, we can deduce $\int_0^Tz_sdB_s=\int_0^T\tilde{z}_sdB_s$. It\^{o} isometry further gives $z_t=\tilde{z}_t,$ $dt\times dP$-a.e., which implies that for each $t\in[0,T]$, $y_t=\tilde{y}_t$. In Definition 3.1(iv), according to Proposition 2.8, we can get that $\widetilde{{\cal{R}}}_T^g(l^\infty_{n\times d}(0,T))\subset\widetilde{{\cal{R}}}_T^0(l^\infty_{n\times d}(0,T))$ and ${\cal{E}}^g[\xi|{\cal{F}}_t]\in\bigcap_{p>1}{\cal{S}}_n^p$ is given in (2.17) explicitly. It is clear that all the cases of Definition 3.1 are consistent. Definition 3.1(ii) is a definition for quadratic $g$-expectations with unbounded terminal variables, a notable example of which is the entropic risk measures (see \cite{BEl} for a bounded variables case). Although Definition 3.1(iii)(iv) contain the superquadratic $g$-expectations with some unbounded terminal variables, in general, the superquadratic $g$-expectations with domain $L^{\infty}({\cal{F}}_T)$ can not be defined clearly, according to the work of Delbaen et al. \cite{DH}.}\\

For convenience, we use $\overline{\cal{G}}_n$ denote the set of all the generators belonging to Definition 3.1, where $n$ is the dimension of generators. When we consider the $g$-expectations in Definition 3.1(ii)(iii), we always assume that $n=1$. We denote the domain of $g$-expectation by $\widetilde{{\cal{R}}}_T^g({\cal{Z}}^g)$. For example, in Definition 3.1(i), ${\cal{Z}}^g={\cal{H}}_{n\times d}^p$. Since $g$ may belong to different cases of Definition 3.1, we use ${\cal{Z}}^g$ denote the largest set among these different sets "${\cal{Z}}^g$". When $n=1$, we denote ${\cal{E}}_n^g$ by ${\cal{E}}^g.$ For $g$-expectations, we have the following properties, whose financial meanings can be found in \cite{R}, \cite{Xu} and the references therein. Note that the following results on the $g$-expectation in Definition 3.1(i) have been obtained in \cite{Peng97}, \cite{R} and \cite{Xu}, etc. We list them for convenience.
\\\\
\textbf{Proposition 3.3}\ \emph{For $g\in\overline{{\cal{G}}}_n$, $\xi,\eta\in\widetilde{{\cal{R}}}_T^g({\cal{Z}}^g)$ and $t\in[0,T]$, we have}

\emph{(i) Monotonicity: ${\cal{E}}^g[\xi|{\cal{F}}_t]\geq{\cal{E}}^g[\eta|{\cal{F}}_t],$ if $n=1$ and $\xi\geq\eta;$}

\emph{(ii) Constant preservation: ${\cal{E}}_n^g[\xi|{\cal{F}}_t]=\xi,$ if $\xi\in\widetilde{{\cal{R}}}_t^g({\cal{Z}}^g);$}

\emph{(iii) Consistency: ${\cal{E}}_n^g[{\cal{E}}_n^g[\xi|{\cal{F}}_r]|{\cal{F}}_t]={\cal{E}}_n^g[\xi|{\cal{F}}_t],$ if $0\leq t\leq r\leq T;$}

\emph{(iv) "0-1 Law": For the $g$-expectations in Definition 3.1(i)-(iii), if $A\in{\mathcal {F}}_t$ such that $1_A\xi\in\widetilde{{\cal{R}}}_T^g({\cal{Z}}^g),$ then} $${\cal{E}}_n^g[1_A\xi|{\cal{F}}_t]=1_A{\cal{E}}_n^g[\xi|{\cal{F}}_t];$$

\emph{(v) Translation invariance: If $g\in{\cal{G}}_n^z$, then for $\zeta\in L_n^0({\mathcal {F}}_t)$ such that $\xi+\zeta\in\widetilde{{\cal{R}}}_T^g({\cal{Z}}^g)$,}
$${\cal{E}}_n^g[\xi+\zeta|{\cal{F}}_t]={\cal{E}}_n^g[\xi|{\cal{F}}_t]+\zeta;$$

\emph{(vi) Concavity: If $n=1$, $g(s,y,z)$ is concave in $(y,z)$, a.e., and for each $\lambda\in[0,1]$, $\lambda\xi+(1-\lambda)\eta\in\widetilde{{\cal{R}}}_T^g({\cal{Z}}^g)$, then for each $\lambda\in[0,1]$,} $${\cal{E}}^g[\lambda\xi+(1-\lambda)\eta|{\cal{F}}_t]\geq\lambda{\cal{E}}^g[\xi|{\cal{F}}_t]+(1-\lambda){\cal{E}}^g[\eta|{\cal{F}}_t];$$

\emph{(vii) Superadditivity: If $n=1$, $g(s,y,z)$ is superadditive in $(y,z)$, a.e., and $\xi+\eta\in\widetilde{{\cal{R}}}_T^g({\cal{Z}}^g)$, then} $${\cal{E}}^g[\xi+\eta|{\cal{F}}_t]\geq{\cal{E}}^g[\xi|{\cal{F}}_t]+{\cal{E}}^g[\eta|{\cal{F}}_t];$$

\emph{(viii) Positive homogeneity: If for each $k>0$ and $(y,z)\in{\mathbf{R}}\times{\mathbf{R}}^{n\times d}$, $g(s,ky,kz)=kg(s,y,z)$, a.e., then for each $k>0$, } $${\cal{E}}_n^g[k\xi|{\cal{F}}_t]= k{\cal{E}}_n^g[\xi|{\cal{F}}_t];$$

\emph{(ix) Independent increments: If $g\in\widehat{{\cal{G}}}_n\cap{\cal{G}}_n^z$, then for each $0\leq t\leq s\leq T$ and $z\in{\mathbf{R}}^{n\times d}$, } $${\cal{E}}_n^g[z(B_s-B_t)|{\cal{F}}_t]={\cal{E}}_n^g[z(B_s-B_t)];$$

\emph{(x) Representation: If $g\in\widehat{{\cal{G}}}_n\cap{\cal{G}}_n^z$, then for each $z\in{\mathbf{R}}^{n\times d}$, } $$g(t,z)=\lim_{\varepsilon\rightarrow0}\frac{1}{\varepsilon}{\cal{E}}_n^g[z(B_{t+\varepsilon}-B_t)],\ \ \textmd{a.e.,}\ t\in[0,T).$$
\emph{Proof.} By Theorem 2.7(i), we get (i). Similar to the arguments in the classical $g$-expectations theory (see \cite{Peng97}), by the uniqueness of solution and the fact that $g(\cdot,\cdot,0)=0$, we can get (ii)-(v).

We prove (vi). By \cite[Proposition 3.5]{EPQ}, we get that the $g$-expectation in Definition 3.1(i) satisfies (vi). Now, we consider the $g$-expectation in Definition 3.1(ii). Since for each $\lambda\in[0,1]$, $\lambda\xi+(1-\lambda)\eta\in\widetilde{{\cal{R}}}_T^g({\cal{H}}^{BMO})$, by Proposition 2.12, we can get that the BSDE($g,T,\lambda\xi+(1-\lambda)\eta$) has a unique $(Y_t,Z_t)\in{\cal{S}} \times{\cal{H}}^{BMO}$ such that $Y_t={\cal{E}}^g[\lambda\xi+(1-\lambda)\eta|{\cal{F}}_t]$. Then by Theorem 2.7(ii) and the arguments as in \cite[Proposition 3.5]{EPQ}, we get that the $g$-expectation in Definition 3.1(ii) satisfies (vi). Similarly, by Theorem 2.7(ii), we can also prove that the $g$-expectations in Definition 3.1(iii)(iv) satisfy (vi).

By a similar argument as in the proof of (vi), we can get (vii). By a similar argument as in the proof of \cite[Proposition 8]{R}, we can get (viii). If $g\in\widehat{{\cal{G}}}_n\cap{\cal{G}}_n^z$, then by Proposition 2.8 and the fact that $g(\cdot,0)=0$, for each $0\leq t\leq s\leq T$ and $z\in{\mathbf{R}}^{n\times d}$,
$${\cal{E}}_n^g[z(B_s-B_t)|{\cal{F}}_t]=\int_t^sg(r,z)dr={\cal{E}}_n^g[z(B_s-B_t)].$$
Thus, (ix) holds. By (ix) and (2.18), we get (x).\ \ $\Box$\\\\
\textbf{Remark 3.4} \emph{The $g$-expectation in Definition 3.1(iv) does not satisfy "0-1" law, except $A=\Omega$ or $A=\emptyset,$ since the domain of such $g$-expectation belongs to $\widetilde{{\cal{R}}}_T^0(l^\infty_{n\times d}(0,T))$ (see Remark 3.2). The assumption: "$\lambda\xi+(1-\lambda)\eta\in\widetilde{{\cal{R}}}_T^g({\cal{Z}}^g)$" in Proposition 3.3(vi) (resp. "$\xi+\eta\in\widetilde{{\cal{R}}}_T^g({\cal{Z}}^g)$" in Proposition 3.3(vii)) is made for the $g$-expectations in Definition 3.1(ii)(iii). }\\

Using Proposition 3.3, we can get the following equivalent results on the comparison, concavity, superadditivity and positive homogeneity for $g$-expectations.\\\\
\textbf{Proposition 3.5} \emph{Let $g, f\in\overline{{\cal{G}}}\cap\widehat{{\cal{G}}}\cap{\cal{G}}^z$ such that ${\cal{Z}}^f\subset{\cal{Z}}^g$ and $\widetilde{{\cal{R}}}_{T}^f({\cal{Z}}^f)\subset\widetilde{{\cal{R}}}_{T}^g({\cal{Z}}^g)$. Then}

\emph{(i) For each $\xi\in\widetilde{{\cal{R}}}_{T}^f({\cal{Z}}^f)$, we have ${\cal{E}}^g[\xi]\geq{\cal{E}}^f[\xi]$ if and only if for each} $z\in{\bf{R}}^d$, $$g(t,z)\geq f(t,z),\ \ \ a.e.,\ t\in[0,T).$$

\emph{(ii) For each $\xi,\eta\in\widetilde{{\cal{R}}}_{T}^g({\cal{Z}}^g)$ such that for each $\lambda\in[0,1]$ $\lambda\xi+(1-\lambda)\eta\in\widetilde{{\cal{R}}}_T^g({\cal{Z}}^g)$, we have ${\cal{E}}^g[\lambda\xi+(1-\lambda)\eta]\geq\lambda{\cal{E}}^g[\xi]+(1-\lambda){\cal{E}}^g[\eta]$
if and only if for each $z_1,z_2\in{\bf{R}}^d$, }
$$g(t,\lambda z_1+(1-\lambda)z_2)\geq \lambda g(t,z_1)+(1-\lambda)g(t,z_2),\ \ \ a.e.,\ t\in[0,T).$$

\emph{(iii) For each $\xi,\eta\in\widetilde{{\cal{R}}}_{T}^g({\cal{Z}}^g)$ such that $\xi+\eta\in\widetilde{{\cal{R}}}_T^g({\cal{Z}}^g)$, we have ${\cal{E}}^g[\xi+\eta]\geq{\cal{E}}^g[\xi]+{\cal{E}}^g[\eta]$ if and only if for each $z_1,z_2\in{\bf{R}}^d$, }
$$g(t,z_1+z_2)\geq g(t,z_1)+g(t,z_2),\ \ \ a.e.,\ t\in[0,T).$$

\emph{(iv) For each $\xi\in\widetilde{{\cal{R}}}_{T}^g({\cal{Z}}^g)$ and $k>0$, we have ${\cal{E}}^g[k\xi]=k{\cal{E}}^g[\xi]$ if and only if for each $k>0$ and $z\in{\bf{R}}^d$,}$$g(t,kz)=kg(t,z),\ \ \ a.e.,\ t\in[0,T).$$
\emph{Proof.} We prove (i). Clearly, for each $(t,z)\in[0,T)\times{\bf{R}}^d$ and $\varepsilon\in[0,T-t]$, $z(B_{t+\varepsilon}-B_t)\in\widetilde{\cal{R}}_T^f({\cal{Z}}^f).$ If for each $(t,z)\in[0,T)\times{\bf{R}}^d$ and $\varepsilon\in[0,T-t]$, we have $${\cal{E}}^g[z(B_{t+\varepsilon}-B_t)]\geq{\cal{E}}^f[z(B_{t+\varepsilon}-B_t)],$$ then by Proposition 3.3(x), we get that for each $z\in{\bf{R}}^d$, $g(t,z)\geq f(t,z)$, $a.e.,$ $t\in[0,T).$

Conversely, for $\xi\in\widetilde{{\cal{R}}}_{T}^f({\cal{Z}}^f)$, since ${\cal{Z}}^f\subset{\cal{Z}}^g$ and $\widetilde{{\cal{R}}}_{T}^f({\cal{Z}}^f)\subset\widetilde{{\cal{R}}}_{T}^g({\cal{Z}}^g)$, by Theorem 2.7(ii), we can get ${\cal{E}}^g[\xi]\geq{\cal{E}}^f[\xi]$.

Similar to the proof of (i), the "right" parts of (ii)(iii)(iv) can be obtained from Proposition 3.3(x). The "left" parts of (ii)(iii)(iv) are due to Proposition 3.3(vi)(vii)(viii), respectively.\ \ $\Box$
\section{${\cal{F}}$-consistent nonlinear expectations}
The definition of filtration-consistent condition expectation (${\cal{F}}$-expectation) was initiated in \cite{CH} in the $1$-dimensional case. We introduce a definition of $n$-dimensional ${\cal{F}}$-expectations. \\\\
\textbf{Definition 4.1} Define a system of operators:
$${\cal{E}}_n[\cdot|{\cal{F}}_t]:\ L_n^0({\cal{F}}_T)\longmapsto L_n^0({\cal{F}}_t), \ t\in[0,T].$$
The operator ${\cal{E}}_n[\cdot|{\cal{F}}_t]$ is called an elementary ${\cal{F}}$-expectation, if it satisfies

(F1) Consistency: ${\cal{E}}_n[{\cal{E}}_n[\xi|{\cal{F}}_t]|{\cal{F}}_s]={\cal{E}}_n[\xi|{\cal{F}}_s],$ if $0\leq s\leq t\leq T;$

(F2) Constant preservation: ${\cal{E}}_n[\xi|{\cal{F}}_t]=\xi,$ if $\xi\in L_n^0({\mathcal {F}}_t).$

An elementary ${\cal{F}}$-expectation ${\cal{E}}_n[\cdot|{\cal{F}}_t]$ is called an ${\cal{F}}$-expectation, if it further satisfies

(F3) Monotonicity: ${\cal{E}}_n[\xi|{\cal{F}}_t]\geq{\cal{E}}_n[\eta|{\cal{F}}_t],$ if $\xi\geq \eta;$

(F4) "0-1 Law": ${\cal{E}}_n[1_A\xi|{\cal{F}}_t]=1_A{\cal{E}}_n[\xi|{\cal{F}}_t],$ if $A\in{\mathcal {F}}_t.$\\

${\cal{F}}$-expectations can be considered as a dynamic counterpart of the risk measures introduced in \cite{Ar} and \cite{FS}, and the dynamic utilities in \cite{DP}. We can also define the concavity, superadditivity and positive homogeneity of ${\cal{F}}$-expectations. From Proposition 3.3, it is clear that on a certain terminal variables set, a $g$-expectation is an elementary ${\cal{F}}$-expectation. However, in general, a $g$-expectation is not an ${\cal{F}}$-expectation in the $n$-dimensional case, since the $n$-dimensional comparison theorem needs some special restrictions on $g$ (see \cite{HP}), Note that we denote ${\cal{E}}_n$ by ${\cal{E}},$ when $n=1$.

Now, we focus on the representation problem described in the Introduction: can we find a function $g\in\overline{{\cal{G}}}_n$ such that
$${\cal{E}}_n[\xi|{\cal{F}}_t]={\cal{E}}_n^g[\xi|{\cal{F}}_t], \ \ \forall t\in[0,T],\eqno(4.1)$$
holds true for $\xi\in\widetilde{{\cal{R}}}_T^g({\cal{Z}}^g)?$

To study the Problem (4.1), we consider the following assumptions:
\begin{itemize}
  \item (A$_{tra}$) (Translation invariance) For each $\xi\in L_n^0({\mathcal {F}}_T)$, $t\in[0,T]$ and $\eta\in L_n^0({\mathcal {F}}_t)$,
  $${\cal{E}}_n[\xi+\eta|{\cal{F}}_t]={\cal{E}}_n[\xi|{\cal{F}}_t]+\eta.$$

  \item (A$_{ind}$) (Independent increments) For each $z\in\textbf{R}^{n\times d}$ and $0\leq t\leq s\leq T,$
  $${\cal{E}}_n[z(B_s-B_t)|{\cal{F}}_t]={\cal{E}}_n[z(B_s-B_t)].$$

  \item (A$_{abs}$) (Absolutely continuous) The function $$G(t,z):={\cal{E}}_n[zB_t],\ \ (t,z)\in[0,T]\times{\mathbf{R}}^{n\times d}$$ is measurable and for each $z\in{\mathbf{R}}^{n\times d}$, $G(t,z)$ is absolutely continuous in $t$.
\end{itemize}

Note that when $n=1$, (A$_{tra}$) has been used in \cite{Co}, \cite{CH}, \cite{Yao} and \cite{Ro}, and (A$_{ind}$) has been used in \cite{Yao}. When (A$_{tra}$) and (A$_{ind}$) hold, (A$_{abs}$) is actually a sufficient and necessary condition for Problem (4.1) (see Theorem 4.3 and Proposition 4.4).

Let ${\cal{Q}}_{n\times d}$ be the set of all the simple step processes defined as
$$\eta_t=\sum^l_{i=0}\sum^{m_i}_{j=1}z_{ij}1_{[t_i,t_{i+1})\times A_{ij}}(t,\omega),\ \ (t,\omega)\in[0,T]\times\Omega,\eqno(4.2)$$
where $0=t_0<t_1<\cdots<t_l<t_{l+1}=T$, $z_{ij}\in{\bf{R}}^{n\times d}$ and $\{A_{ij}\}_{j=1}^{m_i}$ is a partition of $\Omega$ such that for each $0\leq i\leq l$ and $0\leq j\leq m_i$, $A_{ij}\in{\cal{F}}_{t_i}$. It follows that ${\cal{Q}}_{n\times d}$ is dense in ${\cal{H}}_{n\times d}^p,\ p>1$ (This result can be obtained from the proof of \cite[Proposition 2.1]{PR} and a classical approximation for random variables from simple random variables). We further make the following convergence assumptions:

\begin{itemize}
  \item (A$_{L^p}$) For $\xi_m, \xi\in L_n^p({\cal{F}}_T),\ p>1,\ m\geq1,$ if $\xi_m\rightarrow \xi$ in $L_n^p({\cal{F}}_T)$, as $m\rightarrow\infty$, then for each $t\in[0,T]$, ${\cal{E}}_n[\xi_m|{\cal{F}}_t]\rightarrow{\cal{E}}_n[\xi|{\cal{F}}_t]$ in $L_n^p({\cal{F}}_t)$, as $m\rightarrow\infty$.
  \item (A$_{L^\infty}$) For $\xi_m, \xi\in L^\infty({\cal{F}}_T)$ with $\|\xi_m\|_\infty<k$, $m\geq1$ for some $k>0$, if $\xi_m\rightarrow \xi$, $a.s.$, as $m\rightarrow\infty$, then for each $t\in[0,T]$ and $p>1$, ${\cal{E}}[\xi_m|{\cal{F}}_t]\rightarrow{\cal{E}}[\xi|{\cal{F}}_t]$ in $L^p({\cal{F}}_t)$, as $m\rightarrow\infty$.
  \item (A$_{con}$) If $g\in{\cal{L}}^\rho\cap\widehat{{\cal{G}}}\cap{\cal{G}}^z\cap{\cal{G}}^0$ such that (4.1) holds true for $\widetilde{{\cal{R}}}_T^g({\cal{Q}})$, then $g$ is a function such that (4.1) holds true for $\widetilde{{\cal{R}}}_T^g({\cal{Z}}^g)$.
\end{itemize}
\textbf{Remark 4.2} \emph{We make some remarks concerning the convergence assumptions above. The $g$-expectation in Definition 3.1(i) satisfies (A$_{L^p}$) (see \cite[Theorem 5.10]{PR}) and in Definition 3.1(iii) satisfies (A$_{L^\infty}$) (see \cite[Proposition 2.3]{BE}). We will explain that (A$_{con}$) is a natural condition for the $g$-expectations in Definition 3.1(i)(ii)(iii), when $g\in{\cal{G}}^z$. We consider the $g$-expectation in Definition 3.1(ii)(iii). For $g\in{\cal{L}}^\rho\cap\widehat{{\cal{G}}}\cap{\cal{G}}^z\cap{\cal{G}}^0$ and $z\in{\mathbf{R}}^{n\times d}$, there exists $\xi\in\widetilde{{\cal{R}}}_T^g({\cal{Z}}^g)$ such that
$${\cal{E}}^g[\xi|{\cal{F}}_t]={\cal{E}}^g[\xi]-\int_0^tg(s,z)ds+\int_0^tzdB_s,\ \ t\in[0,T].$$
Assume that there exists $\tilde{g}\in{\cal{L}}^\rho\cap\widehat{{\cal{G}}}\cap{\cal{G}}^z\cap{\cal{G}}^0$ such that for each $t\in[0,T]$, ${\cal{E}}^g[\xi|{\cal{F}}_t]={\cal{E}}^{\tilde{g}}[\xi|{\cal{F}}_t].$ Since there exists $Z_s\in{\cal{Z}}^g,$ such that $${\cal{E}}^{\tilde{g}}[\xi|{\cal{F}}_t]={\cal{E}}^{\tilde{g}}[\xi]-\int_0^t{\tilde{g}}(s,Z_s)ds+\int_0^tZ_sdB_s,\ \ t\in[0,T],$$
by comparing the martingale parts of the two equations above, we get that $Z_s=z,\ dt\times dP$-$a.e.,$ and hence for each $z\in{\mathbf{R}}^{n\times d}$, $g(t,z)=\tilde{g}(t,z),\ dt$-$a.e.$
Then it is clear that for each $\xi\in\widetilde{{\cal{R}}}_T^g({\cal{Z}}^g)$, we have ${\cal{E}}^g[\xi|{\cal{F}}_t]={\cal{E}}^{\tilde{g}}[\xi|{\cal{F}}_t],$
which implies that ${\cal{E}}^{\tilde{g}}[\xi|{\cal{F}}_t]$ satisfies (A$_{con}$) on $\widetilde{{\cal{R}}}_T^g({\cal{Z}}^g).$ Similarly, we can get that the $g$-expectations in Definition 3.1(i) also satisfy (A$_{con}$).} \\

We have the following representation theorem under absolutely continuous condition.\\\\
\textbf{Theorem 4.3} \emph{Let the elementary ${\cal{F}}$-expectation ${\cal{E}}_n$ satisfy (A$_{ind}$), (A$_{tra}$) and (A$_{abs}$). Then there exists a unique function $g\in\widehat{{\cal{G}}}_n\cap{\cal{G}}_n^z\cap{\cal{G}}_n^0$ such that (4.1) holds true for ${\cal{R}}_n$. Moreover if we assume that ${\cal{E}}_n$ also satisfies (F4), then}

\emph{(i) (4.1) holds true for $L_n^p({\cal{F}}_T)$, when $g\in{{\cal{L}}}_n^\nu$, $p>1$ and (A$_{L^p}$) holds;}

\emph{(ii) (4.1) holds true for $\widetilde{{\cal{R}}}_T^g({\cal{H}}^\infty)$, in particular for $\overline{{\cal{R}}}$, when $g\in{\cal{L}}^\rho$, $p>1$ and (A$_{L^p}$) holds;}

\emph{(iii) (4.1) holds true for $\widetilde{{\cal{R}}}_T^g({\cal{Z}}^g)$, in particular for $L^{\infty}({\cal{F}}_T)$ in the case that $\rho(|\cdot|)=\nu(|\cdot|+1)$, when $g\in{\cal{L}}^\rho$ and (A$_{con}$) holds.}
\\\\
\emph{Proof.} \textbf{Step 1.} Since ${\cal{E}}_n$ satisfies (A$_{abs}$), there exists a function $g(t,z)$ defined on $[0,T]\times{\mathbf{R}}^{n\times d}$ such that for each $z\in{\mathbf{R}}^{n\times d},$
$$g(t,z)=\partial_tG(t,z)=\partial_t{\cal{E}}_n[zB_t],\ \ \ a.e.\ t\in(0,T).$$
It is easy to get that $g\in\widehat{{\cal{G}}}_n\cap{\cal{G}}_n^z\cap{\cal{G}}_n^0$. For each $z\in{\mathbf{R}}^{n\times d}$ and $0\leq t\leq s\leq T$, by (A$_{abs}$), (A$_{tra}$), (F1) and (A$_{ind}$), we have
\begin{eqnarray*}
\ \ \ \ \ \ \ \ \ \ \ \ \ \ \ \ \ \ \ \ \ \ \ \ \int_t^sg(r,z)dr&=&{\cal{E}}_n[zB_s]-{\cal{E}}_n[zB_t]\\&=&{\cal{E}}_n[zB_s-{\cal{E}}_n[zB_t]]\\
&=&{\cal{E}}_n[{\cal{E}}_n[z(B_s-B_t)+zB_t-{\cal{E}}_n[zB_t]|{\cal{F}}_t]]\\
&=&{\cal{E}}_n[{\cal{E}}_n[z(B_s-B_t)|{\cal{F}}_t]+zB_t-{\cal{E}}_n[zB_t]]\\
&=&{\cal{E}}_n[z(B_s-B_t)]+{\cal{E}}_n[zB_t-{\cal{E}}_n[zB_t]]\\
&=&{\cal{E}}_n[z(B_s-B_t)].\ \ \ \ \ \ \ \ \ \ \ \ \ \ \ \ \ \ \ \ \ \ \ \ \ \ \ \ \ \ \ \ \ \ \ \ \ \ \ \ \ \ \ \ \ \ \ \ \ \ \ (4.4)
\end{eqnarray*}
Then by (A$_{tra}$), (A$_{ind}$) and (4.4), we get that for each $z\in{\mathbf{R}}^{n\times d}$ and $0\leq t\leq s\leq T$,
\begin{eqnarray*}
\ \ \ \ \ \ \ \ \ \ \ \ {\cal{E}}_n\left[-\int_t^sg(r,z)dr+\int_t^szdB_r|{\cal{F}}_t\right]&=&-\int_t^sg(s,z)dr+{\cal{E}}_n[z(B_s-B_t)|{\cal{F}}_t]\\
&=&-\int_t^sg(s,z)dr+{\cal{E}}_n[z(B_s-B_t)]\\
&=&0.\ \ \ \ \ \ \ \ \ \ \ \ \ \ \ \ \ \ \ \ \ \ \ \ \ \ \ \ \ \ \ \ \ \ \ \ \ \ \ \ \ \ \ \ \ \ \ \ \ \ \ \ (4.5)
\end{eqnarray*}
For $\xi\in{\cal{R}}_n$, there exist $(y,z)\in{\mathbf{R}}^n\times{\mathbf{R}}^{n\times d}$ and $0\leq u\leq v\leq T$ such that $\xi=y+z(B_v-B_u).$ Since $g\in{\cal{G}}_n^0$, we get $\int_0^T|g(r,z1_{[u,v]}(r))|dr<\infty.$ Then by Proposition 2.8 and Definition 3.1(iv), we have
$${\cal{E}}_n^g[\xi|{\cal{F}}_t]=\xi+\int_t^Tg(r,z1_{[u,v]}(r))dr-\int_t^Tz1_{[u,v]}(r)dB_r,\ \ t\in[0,T].$$
Then by (A$_{tra}$), (F1), (F2), (4.5) and the fact $g\in{\cal{G}}_n^0$, we can deduce that
\begin{eqnarray*}
 \ \ \ \ \ \ \ \ {\cal{E}}_n[\xi|{\cal{F}}_t]-{\cal{E}}_n^g[\xi|{\cal{F}}_t]
&=&{\cal{E}}_n\left[-\int_t^Tg(r,z1_{[u,v]}(r))dr+\int_t^Tz1_{[u,v]}(r)dB_r|{\cal{F}}_t\right]\\
&=&\left\{
       \begin{array}{ll}
         {\cal{E}}_n\left[{\cal{E}}_n\left[-\int_u^vg(r,z)dr+\int_u^vzdB_r|{\cal{F}}_u\right]|{\cal{F}}_t\right], \ \ &t\in[0,u];\\
         {\cal{E}}_n\left[-\int_t^vg(r,z)dr+\int_t^vzdB_r|{\cal{F}}_t\right], \ \ \ &t\in(u,v);\\
0,\ \ \ &t\in[v,T]
       \end{array}
     \right.\\
&=&0,\ \ \ t\in[0,T]. \ \ \ \ \ \ \ \ \ \ \ \ \ \ \ \ \ \ \ \ \ \ \ \ \ \ \ \ \ \ \ \ \ \ \ \ \ \ \ \
 \ \ \ \ \ \ \ \ \ \ \ \ \ \ \ \ \ \ \ \ (4.6)
\end{eqnarray*}
which implies that $g$ is a function in $\widehat{{\cal{G}}}_n\cap{\cal{G}}_n^z\cap{\cal{G}}_n^0$ such that (4.1) holds for ${\cal{R}}_n$.

If there exists another function $\bar{g}\in\widehat{{\cal{G}}}_n\cap{\cal{G}}_n^z\cap{\cal{G}}_n^0$ such that for each $\xi\in{\cal{R}}_n$,
$${\cal{E}}_n[\xi|{\cal{F}}_t]={\cal{E}}_n^{\bar{g}}[\xi|{\cal{F}}_t],\ \ \forall t\in[0,T],$$
then by (4.6),  for each $z\in{\mathbf{R}}^{n\times d}$, we have,
$${\cal{E}}_n^g[z(B_{t+\varepsilon}-B_t)|{\cal{F}}_t]={\cal{E}}_n^{\bar{g}}[z(B_{t+\varepsilon}-B_t)|{\cal{F}}_t],\ \ \forall t\in[0,T].\eqno(4.7)$$
By applying Proposition 3.1(ix)(x) to (4.7), we get that for each $z\in\textbf{R}^d$,
$$g(t,z)=\bar{g}(t,z),\ \ \ a.e.,\ t\in[0,T).\eqno(4.8)$$
Thus, $g$ is a unique function in $\widehat{{\cal{G}}}_n\cap{\cal{G}}_n^z\cap{\cal{G}}_n^0$ such that (4.1) holds true for ${\cal{R}}_n$.

\textbf{Step 2.} We assume that (F4) holds true. For each $0\leq r<t\leq T$, let $\eta:=\sum_{j=1}^m1_{A_j}z_j$, where $z_j\in{\bf{R}}^{n\times d}$ and $\{A_j\}_{j=1}^m$ is a partition of $\Omega$ with $A_j\in{\cal{F}}_r, 1\leq j\leq m$. Then by (F4) and (4.5), we have
\begin{eqnarray*}
\ \ \ \ \ {\cal{E}}_n\left[-\int_r^tg(s,\eta)ds+\int_r^t\eta dB_s|{\cal{F}}_r\right]&=&{\cal{E}}_n\left[\sum_{j=1}^m 1_{A_j}(-\int_r^tg(s,z_j)ds+\int_r^tz_j dB_s)|{\cal{F}}_r\right]\\&=&\sum_{j=1}^m 1_{A_j}{\cal{E}}_n\left[\sum_{j=1}^m 1_{A_j}(-\int_r^tg(s,z_j)ds+\int_r^tz_j dB_s)|{\cal{F}}_r\right]\\&=&\sum_{j=1}^m{\cal{E}}_n\left[1_{A_j}(-\int_r^tg(s,z_j)ds+\int_r^tz_j dB_s)|{\cal{F}}_r\right]\\&=&\sum_{j=1}^m 1_{A_j}{\cal{E}}_n\left[-\int_r^tg(s,z_j)ds+\int_r^tz_j dB_s|{\cal{F}}_r\right]\\
&=&0. \ \ \ \ \ \ \ \ \ \ \ \ \ \ \ \ \ \ \ \ \ \ \ \ \ \ \ \ \ \ \ \ \ \ \ \ \ \ \ \ \ \ \ \ \ \ \ \ \ \ \ \ \ \ \ \ \ \ \ (4.9)
\end{eqnarray*}
For any simple step process $\eta_t\in{\cal{Q}}_{n\times d}$ (see (4.2)) and $t\in[0,T]$, if $t\in[t_i,t_{i+1})$, then we set $r_i=t$ and $r_k=t_k,$ $i<k\leq l+1$. Then by (A$_{tra}$) and (4.9), we have
\begin{eqnarray*}
\ \ \ \ \ \ \ \ \ \ \ \ \ \ \ &&{\cal{E}}_n\left[-\int_t^Tg(s,\eta_s)ds+\int_t^T\eta_sdB_s|{\cal{F}}_t\right]\\
&=&{\cal{E}}_n\left[\sum^l_{k=i}\left(-\int_{r_k}^{r_{k+1}}g(s,\sum_{j=1}^{m_k}z_{kj}1_{A_{kj}})ds
+\int_{r_k}^{r_{k+1}}\sum_{j=1}^{m_k}z_{kj}1_{A_{kj}}dB_s\right)|{\cal{F}}_t\right]\\
&=&{\cal{E}}_n\left[{\cal{E}}_n\left[\left(-\int_{r_l}^{T}g(s,\sum_{j=1}^{m_l}z_{lj}1_{A_{lj}})ds
+\int_{r_l}^{T}\sum_{j=1}^{m_l}z_{lj}1_{A_{lj}}dB_s\right)|{\cal{F}}_{r_l}\right]\right.\\
&&\left.+\sum^{l-1}_{k=i}\left(-\int_{r_k}^{r_{k+1}}g(s,\sum_{j=1}^{m_k}z_{kj}1_{A_{kj}})ds
+\int_{r_k}^{r_{k+1}}\sum_{j=1}^{m_k}z_{kj}1_{A_{kj}}dB_s\right)|{\cal{F}}_t\right]
\\&=&{\cal{E}}_n\left[\sum^{l-1}_{k=i}\left(-\int_{r_k}^{r_{k+1}}g(s,\sum_{j=1}^{m_k}z_{kj}1_{A_{kj}})ds
+\int_{r_k}^{r_{k+1}}\sum_{j=1}^{m_i}z_{kj}1_{A_{kj}}dB_s\right)|{\cal{F}}_t\right]\\
&\cdots&\\&=&{\cal{E}}_n\left[-\int_{t}^{t_{i+1}}g(s,\sum_{j=1}^{m_i}z_{ij}1_{A_{ij}})ds
+\int_{t}^{t_{i+1}}\sum_{j=1}^{m_i}z_{ij}1_{A_{ij}}dB_s|{\cal{F}}_t\right]\\
&=&0.\ \ \ \ \ \ \ \ \ \ \ \ \ \ \ \ \ \ \ \ \ \ \ \ \ \ \ \ \ \ \ \ \ \ \ \ \ \ \ \ \ \ \ \ \ \ \ \ \ \ \ \ \ \ \ \ \ \ \ \ \ \ \ \ \ \ \ \ \ \ \ \ \ \ \ \ \ \ \ \ \ \ \ \ \ \ \ \ \ \ (4.10)
\end{eqnarray*}

\textbf{Step 3.} We prove (i). Let $p>1$. If $g\in{\cal{L}}_n^\nu$, then by Definition 3.1(i), for each $\xi\in L_n^p({\cal{F}}_T)$, we have
$${\cal{E}}_n^g[\xi|{\cal{F}}_t]=\xi+\int_t^Tg(s,Z_s)ds-\int_t^TZ_sdB_s,\ \ t\in[0,T],\eqno(4.11)$$
with $Z_t\in{\cal{H}}_{n\times d}^p$. Since ${\cal{Q}}_{n\times d}$ is dense in ${\cal{H}}_{n\times d}^p$, there exists a series of simple step processes $\{\eta^l_t\}_{l\geq1}\in{\cal{Q}}_{n\times d}$ such that $\eta_t^l\rightarrow Z_t$ in ${\cal{H}}_{n\times d}^p$ as $l\rightarrow\infty$. By the Lipschitz condition, we can get that for each $t\in[0,T]$,
$$-\int_t^Tg(s,\eta_s^l)ds+\int_t^T\eta_s^ldB_s\longrightarrow-\int_t^Tg(s,Z_s)ds+\int_t^TZ_sdB_s$$
in $L_n^p$, as $l\rightarrow\infty$. Then by (A$_{L^p}$) and (4.10), we have
$${\cal{E}}_n\left[-\int_t^Tg(s,Z_s)ds+\int_t^TZ_sdB_s|{\cal{F}}_t\right]=0,\ \ \forall t\in[0,T], \eqno(4.12)$$
which together with (4.11) and (A$_{tra}$) implies (4.1) holds true for $\xi.$ Thus, (i) holds true.

We prove (ii). Let $p>1$. For $\xi\in\widetilde{{\cal{R}}}_T^g({\cal{H}}^\infty)$, there exists a constant $k>0$ such that $\xi\in\widetilde{{\cal{R}}}_T^g({\cal{H}}^{\infty,k})$. If $g\in{\cal{L}}^\rho$, then by Definition 3.1(iii), we have
$${\cal{E}}^{g}[\xi|{\cal{F}}_t]=\xi+\int_t^Tg(s,Z_s)ds-\int_t^TZ_sdB_s,\ \ t\in[0,T],$$
with $Z_t\in{\cal{H}}^{\infty,k}$. Moreover, we can find a series of simple step processes $\{\eta^l_t\}_{l\geq1}\in{\cal{Q}}$ such that for each $l\geq1$, $\|\eta^l\|_\infty\leq k$ and $\eta_t^l\rightarrow Z_t$ in ${\cal{H}}^p$, as $l\rightarrow\infty$ (This can be proved, according to the proof of \cite[Proposition 2.1]{PR} and a classical approximation for bounded variables by simple variables with the same bound). Thus, by the same argument as in (4.12), we can get
$${\cal{E}}\left[-\int_t^Tg(s,Z_s)ds+\int_t^TZ_sdB_s|{\cal{F}}_t\right]=0,\ \ \forall t\in[0,T].$$
By the above two equalities and (A$_{tra}$), we get that (4.1) holds true for $\widetilde{{\cal{R}}}_T^g({\cal{H}}^\infty),$ in particular for $\overline{{\cal{R}}}$ (see Proposition 2.5(iv)).

We prove (iii). If $g\in{\cal{L}}^{\rho}$, then by (4.10) and (A$_{tra}$), we can get that (4.1) holds true for $\widetilde{{\cal{R}}}_T^g({\cal{Q}})$. Moreover, since (A$_{con}$) holds true, we get that $g$ is a function such that (4.1) holds true for $\widetilde{{\cal{R}}}_T^g({\cal{Z}}^g)$, in particular for $L^{\infty}({\cal{F}}_T)$ in the case that $\rho(|\cdot|)=\nu(|\cdot|+1)$ (see Definition 3.1(ii) and Proposition 2.5(ii)).

The proof is complete.\ \ $\Box$\\

The following Proposition 4.4 shows that (A$_{abs}$) is a necessary condition for Problem (4.1) under (A$_{int}$) and (A$_{tra}$).\\\\
\textbf{Proposition 4.4} \emph{Let the elementary ${\cal{F}}$-expectation ${\cal{E}}_n$ satisfy (A$_{ind}$) and (A$_{tra}$). If there exists a function $g\in\widehat{{\cal{G}}}_n\cap{\cal{G}}_n^z\cap{\cal{G}}_n^0$ such that (4.1) holds true for ${\cal{R}}_n$, then ${\cal{E}}_n$ satisfies (A$_{abs}$).}\\\\
\emph{Proof.} Since the elementary ${\cal{F}}$-expectation ${\cal{E}}_n$ satisfy (A$_{ind}$) and (A$_{tra}$), according to the arguments in (4.4), we have for each $z\in{\mathbf{R}}^{n\times d}$,
$${\cal{E}}_n[z(B_T-B_t)]={\cal{E}}_n[zB_T]-{\cal{E}}_n[zB_t],\ \ \forall t\in[0,T].\eqno(4.13)$$
In view of (2.17), for each $g\in\widehat{{\cal{G}}}_n\cap{\cal{G}}_n^z\cap{\cal{G}}_n^0$, we have
$${\cal{E}}_n^g[z(B_T-B_t)]=\int_t^Tg(s,z)ds,\ \ \forall t\in[0,T].\eqno(4.14)$$
Since there exists a function $g\in\widehat{{\cal{G}}}_n\cap{\cal{G}}_n^z\cap{\cal{G}}_n^0$ such that (4.1) holds true for ${\cal{R}}_n$, by (4.13) and (4.14), we get $${\cal{E}}_n[zB_t]=\int_0^tg(s,z)ds,\ \ \forall t\in[0,T],$$ which implies that ${\cal{E}}_n$ satisfies (A$_{abs}$).\ \ $\Box$\\\\
\textbf{Remark 4.5} \emph{In Theorem 4.3, the results are obtained without using the Monotonicity (F3) and the domination method developed by \cite{CH}, and the generator $g$ can be determined without any assumptions on $z$. Theorem 4.3(i) contains a representation theorem for $n$-dimensional $\cal{F}$-expectations. This result can not be fully obtained using the domination method (a special result using domination method was announced in the arXiv version of \cite{Xu}), since the $n$-dimensional comparison theorem for BSDEs is not true in general (see \cite{HP}). Theorem 4.3(iii) contains a representation theorem for quadratic $\cal{F}$-expectations, a related result for bounded terminal variables without assuming a condition like (A$_{con}$), was obtained by \cite{Yao} using the domination method under (A$_{int}$), (A$_{tra}$) and several domination conditions (see \cite[Definition 3.8(1)-(3) and H4]{Yao}).}\\

Now, we study the Problem (4.1) in the $1$-dimensional case under the following domination condition (A$_{\rho(k)}$), $k>0$. Note that for $k>0$, we denote ${\cal{E}}^{g}$ by ${\cal{E}}^{k}$ (resp. denote ${\cal{E}}^{g}$ by ${\cal{E}}^{-k}$),  if $g=k|z|$  (resp. $g=-k|z|$).

\begin{itemize}
  \item (A$_{\rho(k)}$) ($\rho(k)$-domination on ${\cal{K}}$)  For $k>0$, $p>1$ and $\xi, \eta\in{\cal{K}}\subset L^p({\cal{F}}_T)$, we have
\begin{center}
${\cal{E}}[\xi|{\cal{F}}_t]-{\cal{E}}[\eta|{\cal{F}}_t]\leq{\cal{E}}^{\rho(k)}[\xi-\eta|{\cal{F}}_t],\ \ \forall t\in[0,T].$
\end{center}
\end{itemize}

In fact, from \cite[Lemma 4.3 and Lemma 4.4]{CH} and the fact ${\cal{E}}^{\rho(k)}[\xi|{\cal{F}}_t]=-{\cal{E}}^{\rho(k)}[-\xi|{\cal{F}}_t]$, it follows that (A$_{\rho(k)}$) on $L^2({\cal{F}}_T)$ is equivalent to the original domination condition:
$${\cal{E}}[\xi]-{\cal{E}}[\eta]\leq{\cal{E}}^{\rho(k)}[\xi-\eta],\ \ \forall \xi,\eta\in L^2({\cal{F}}_T)$$
plus the strict monotonicity condition and translation invariance condition in \cite{CH} (see \cite{Co} and \cite{Ro} in general filtrations and \cite{Zheng1} for a generalization). We point out that (A$_{\rho(k)}$) is also a sufficient condition for (A$_{tra}$) and (A$_{L^p}$) (see Lemma A.2(ii)(iii) in the Appendix). Under (A$_{\rho(k)}$), we have the following representation theorem, whose proof relies on a special Doob-Meyer decomposition for ${\cal{E}}$-supermartingales given in the Appendix.\\\\
\textbf{Theorem 4.6} \emph{Let the ${\cal{F}}$-expectation ${\cal{E}}$ satisfy (A$_{ind}$) and (A$_{\rho(k)}$) on ${\cal{R}}^k$ for each $k>0$. Then there exists a unique function $g\in\widehat{{\cal{G}}}\cap{\cal{G}}^z\cap{\cal{G}}^0\cap{{\cal{L}}}^{\rho}$ such that (4.1) holds true for ${\cal{R}}$. Moreover, if we assume that ${\cal{E}}$ satisfies (A$_{\rho(k)}$) on $\widetilde{{\cal{R}}}_T^g({\cal{H}}^{\infty,k})$ for each $k>0$, then $g$ is a function such that (4.1) holds true for $\widetilde{{\cal{R}}}_T^g({\cal{H}}^{\infty})$, in particular for $\overline{{\cal{R}}}$, and when $\rho(|\cdot|)=\nu(|\cdot|+1)$ and (A$_{L^\infty}$) holds true, $g$ is a function such that (4.1) holds true for $L^{\infty}({\cal{F}}_T)$.}\\\\
\emph{Proof.} \textbf{Step 1.} For $k>0$ and $z\in\textbf{B}(k)$, by Definition 3.1(i), it is not hard to check that $(-\rho(k)|z|t+zB_t,z)$ is a unique solution to the BSDE$(\rho(k)|z|,T,-\rho(k)|z|T+zB_T)$ with $$-\rho(k)|z|t+zB_t={\cal{E}}^{\rho(k)}[-\rho(k)|z|T+zB_T|{\cal{F}}_t],\ \ \forall t\in[0,T].\eqno(4.15)$$ Since the ${\cal{F}}$-expectation ${\cal{E}}$ satisfies (A$_{\rho(k)}$) on ${\cal{R}}^k$, it follows from (4.15) and Lemma A.2(i) that $-\rho(k)|z|t+zB_t$ is an ${\cal{E}}$-supermartingale (the notion of ${\cal{E}}$-supermartingale is introduced in the Appendix). Then by Lemma A.6, there exists a function $a(t,z)\in l^2_{1\times1}(0,T)$, which is continuous and increasing with $a(0,z)=0$ such that for each $t\in[0,T],$
$${\cal{E}}[-\rho(k)|z|T+zB_T+a(T,z)|{\cal{F}}_t]=-\rho(k)|z|t+zB_t+a(t,z).\eqno(4.16)$$
By Lemma A.4, there exists a pair $(g_k(t,z), Z^z_t)\in{\cal{H}}_{1\times1}^2\times{\cal{H}}_{1\times d}^2$ such that for each $t\in[0,T],$
$$-\rho(k)|z|t+zB_t+a(t,z)=-\rho(k)|z|T+zB_T+a(T,z)+\int_t^Tg_k(r,z)dr-\int_t^TZ_r^zdB_r,\eqno(4.17)$$
and for each $t\in[0,T]$ and $\bar{z}\in \textbf{B}(k)$,
$$|g_k(t,z)|\leq\rho(k)|Z_t^z| \ \ \textrm{and}\ \ |g_k(t,z)-g_k(t,\bar{z})|\leq \rho(k)|Z_t^z-Z_t^{\bar{z}}|,\ \ dt\times dP-a.e.\eqno(4.18)$$
By comparing the martingale parts and the bounded variation parts in (4.17), we get
$$Z_t^z=z\ \ dt\times dP-a.e.\ \ \textmd{and}\ \ \int_0^tg_k(r,z)dr=\rho(k)|z|t-a(t,z),\ \ \forall t\in[0,T].\eqno(4.19)$$
By (4.18), (4.19) and the fact that $a(t,z)$ is deterministic, we can get $g_k(t,z)\in\widehat{{\cal{G}}}\cap{\cal{G}}^z\cap{\cal{G}}^0\cap{{\cal{L}}}^{\rho(k)}$ on $[0,T]\times{\bf{B}}(k).$ By (4.16), (4.17), (4.19) and Lemma A.2(iii), for each $k>0$ and $z\in{\bf{B}}(k)$, we deduce that
$${\cal{E}}\left[-\int_s^tg_k(r,z)ds+\int_s^tzdB_r|{\cal{F}}_s\right]=0,\ \ \forall 0\leq s\leq t\leq T.\eqno(4.20)$$
Then by (4.20), Lemma A.2(iii) and a similar argument as in (4.6), we get for each $k>0$ and $\xi\in{\cal{R}}^k$,
$${\cal{E}}[\xi|{\cal{F}}_t]={\cal{E}}^{g_k}[\xi|{\cal{F}}_t],\ \ \forall t\in[0,T],\eqno(4.21)$$
which implies that for each $k>0$, $r>0$ and $\xi\in{\cal{R}}^k$,
$${\cal{E}}^{g_k}[\xi|{\cal{F}}_t]={\cal{E}}^{g_{k+r}}[\xi|{\cal{F}}_t],\ \ \forall t\in[0,T].$$
Then by a similar argument as in (4.8), we get for each $k>0$, $r>0$ and $z\in{\bf{B}}(k)$,
$$g_k(t,z)=g_{k+r}(t,z),\ \ a.e.\ t\in[0,T).$$
For $(t,z)\in[0,T]\times{\bf{R}}^d$, set
$$g(t,z):=\left\{
    \begin{array}{ll}
      g_1(t,z), &(t,z)\in[0,T]\times{\bf{B}}(1); \\
      g_m(t,z), & (t,z)\in[0,T]\times{\bf{B}}(m)\}\setminus[0,T]\times{\bf{B}}(m-1), \ \ m\geq2.
    \end{array}
  \right.$$
Thus, $g(t,z)$ is well-defined on $[0,T]\times{\bf{R}}^d$. Moreover, for each $k>0$, we have $g(t,z)\in\widehat{{\cal{G}}}\cap{\cal{G}}^z\cap{\cal{G}}^0\cap{{\cal{L}}}^{\rho(k)}$ on $[0,T]\times{\bf{B}}(k),$ which implies that for each $z_i\in{\mathbf{R}}^{\mathit{d}},\ i=1,2,$
  $$|{g}(t,z_{1})-{g}(t,z_{2}) |\leq\rho(|z_1|\vee|z_2|)|z_{\mathrm1}-z_{2}|,\ \ a.e.\ t\in[0,T).$$
Thus, $g(t,z)\in\widehat{{\cal{G}}}\cap{\cal{G}}^z\cap{\cal{G}}^0\cap{\cal{L}}^{\rho}$. For $\xi\in{\cal{R}}$, there exists a constant $k>0$ such that $\xi\in{\cal{R}}^k.$ Then by Proposition 2.8, we get $\xi\in\widetilde{{\cal{R}}}_T^g({\cal{H}}^{\infty,k})=\widetilde{{\cal{R}}}_T^{g_k}({\cal{H}}^{\infty,k}).$ Hence, by Definition 3.1(i)(iii), we have for each $t\in[0,T]$, ${\cal{E}}^g[\xi|{\cal{F}}_t]={\cal{E}}^{g_k}[\xi|{\cal{F}}_t]$, which together with (4.21) implies that for each $\xi\in{\cal{R}}$,
$${\cal{E}}[\xi|{\cal{F}}_t]={\cal{E}}^g[\xi|{\cal{F}}_t],\ \ \forall t\in[0,T].\eqno(4.22)$$

If there exists another function $\hat{g}(t,z)\in\widehat{{\cal{G}}}\cap{\cal{G}}^z\cap{\cal{G}}^0\cap{{\cal{L}}}^{\rho}$ such that for each $\xi\in{\cal{R}}$,
$${\cal{E}}[\xi|{\cal{F}}_t]={\cal{E}}^{\hat{g}}[\xi|{\cal{F}}_t],\ \ \forall t\in[0,T].\eqno(4.23)$$
Then similar to the argument in (4.8), by (4.22) and (4.23), we get that for each $z\in{\bf{R}}^{d},$
$$g(t,z)=\bar{g}(t,z),\ \ a.e.\ t\in[0,T).$$
This means that $g$ is a unique function in $\widehat{{\cal{G}}}\cap{\cal{G}}^z\cap{\cal{G}}^0\cap{{\cal{L}}}^{\rho}$ such that (4.1) holds true for ${\cal{R}}$.

\textbf{Step 2.} We assume that ${\cal{E}}$ satisfies (A$_{\rho(k)}$) on $\widetilde{{\cal{R}}}_T^g({\cal{H}}^{\infty,k})$ for each $k>0$. Since $g\in{\cal{G}}^0,$ we can get that for each $\eta_t\in{\cal{Q}}\cap{\cal{H}}^{\infty,k}$ and $0\leq t\leq s\leq T$, $-\int_t^sg(r,\eta_r)dr+\int_t^s\eta_rdB_r$ belongs to $\widetilde{{\cal{R}}}_T^g({\cal{H}}^{\infty,k}).$ Then by (4.20) and the definition of $g$, Lemma A.2(iii) and the Step 2 in the proof of Theorem 4.3, we can get that for each $\eta_t\in{\cal{Q}}\cap{\cal{H}}^{\infty,k}$,
$${\cal{E}}\left[-\int_t^Tg(s,\eta_s)ds+\int_t^T\eta_sdB_s|{\cal{F}}_t\right]=0,\ \ \forall t\in[0,T]. \eqno(4.24)$$
For $\xi\in\widetilde{{\cal{R}}}_T^g({\cal{H}}^{\infty})$, there exists a constant $k>0$ such that $\xi\in\widetilde{{\cal{R}}}_T^g({\cal{H}}^{\infty,k})$. By Definition 3.1(iii), we have
$${\cal{E}}^g[\xi|{\cal{F}}_t]=\xi+\int_t^Tg(s,Z_s)ds-\int_t^TZ_sdB_s,\ \ t\in[0,T],\eqno(4.25)$$
with $Z_t\in{\cal{H}}^{\infty,k}$. Then similar to the argument in the proof of Theorem 4.3(ii), we can find a series of simple step processes $\{\eta^l_t\}_{l\geq1}\subset{\cal{Q}}$ such that for each $l\geq1,$ $\|\eta_t^l\|_{\infty}\leq k$ and $\eta_t^l\rightarrow Z_t$ in ${\cal{H}}_{n\times d}^2$ as $l\rightarrow\infty$. By the locally Lipschitz condition, we can get that
$$-\int_t^Tg(s,\eta_s^l)ds+\int_t^T\eta_s^ldB_s\longrightarrow-\int_t^Tg(s,Z_s)ds+\int_t^TZ_sdB_s, \eqno(4.26)$$
in $L^2$, as $l\rightarrow\infty$. Since $g\in{\cal{G}}^0,$ we can get that $-\int_t^Tg(s,\eta_s^l)ds+\int_t^T\eta_s^ldB_s$ and $-\int_t^Tg(s,Z_s)ds+\int_t^TZ_sdB_s$ both belong to $\widetilde{{\cal{R}}}_T^g({\cal{H}}^{\infty,k}).$ Moreover, since ${\cal{E}}$ satisfies (A$_{\rho(k)}$) on $\widetilde{{\cal{R}}}_T^g({\cal{H}}^{\infty,k})$ for each $k>0$, by Lemma A.2(ii), we have
\begin{eqnarray*}
&&\left|{\cal{E}}\left[-\int_t^Tg(s,\eta_s^l)ds+\int_t^T\eta_s^ldB_s|{\cal{F}}_t\right]-{\cal{E}}\left[-\int_t^Tg(s,Z_s)ds+\int_t^TZ_sdB_s|{\cal{F}}_t\right]\right|\\
&\leq&{\cal{E}}^{\rho(k)}\left[\left|-\int_t^Tg(s,\eta_s^l)ds+\int_t^T\eta_s^ldB_s-\left(-\int_t^Tg(s,Z_s)ds+\int_t^TZ_sdB_s\right)\right||{\cal{F}}_t\right],\ \ \forall t\in[0,T],
\end{eqnarray*}
which together with (4.26), \cite[Theorem 3.2]{Peng04} and (4.24) implies
$${\cal{E}}\left[-\int_t^Tg(s,Z_s)ds+\int_t^TZ_sdB_s|{\cal{F}}_t\right]=0,\ \ \forall t\in[0,T].$$
This together with (4.25) and Lemma A.2(iii) implies that (4.1) holds for $\xi$. Thus, $g$ is a function such that (4.1) holds true for $\widetilde{{\cal{R}}}_T^g({\cal{H}}^{\infty})$, in particular for $\overline{{\cal{R}}}$.

We further assume that $\rho(|\cdot|)=\nu(|\cdot|+1)$ and (A$_{L^\infty}$) holds. Since $\rho(|\cdot|)=\nu(|\cdot|+1)$, by the proof of \cite[Theorem 2.2]{BE}, for each $\xi\in L^{\infty}({\cal{F}}_T)$, we can find a series $\{\xi_m\}_{m\geq1}\subset\overline{{\cal{R}}}$, which is a Cauchy sequence in $L^2$ such that $\sup_{m\geq1}\|\xi_m\|_\infty\leq\|\xi\|_\infty$ and for each $t\in[0,T]$, ${\cal{E}}^g[\xi_m|{\cal{F}}_t]\rightarrow{\cal{E}}^g[\xi|{\cal{F}}_t]$ in $L^p$ for each $p>1$, as $m\rightarrow \infty.$ This result together with (A$_{L^\infty}$) and the fact that $g$ is a function such that (4.1) holds true for $\overline{{\cal{R}}}$, implies that $g$ is a function such that (4.1) holds true for $L^{\infty}({\cal{F}}_T).$\ \  $\Box$\\\\
\textbf{Corollary 4.7} \emph{Let $k>0$ and the ${\cal{F}}$-expectation ${\cal{E}}$ satisfy (A$_{int}$) and (A$_{\rho(k)}$) on $L^2({\cal{F}}_T).$ Then there exists a unique function $g\in\widehat{{\cal{G}}}\cap{\cal{G}}^z\cap{\cal{G}}^0\cap{\cal{L}}^{\rho(k)}$ such that (4.1) holds true for $L^2({\cal{F}}_T)$.}\\\\
\emph{Proof.} Since ${\cal{R}}\in L^2({\cal{F}}_T)$, by Theorem 4.6, we can get that there exists a unique function $g\in\widehat{{\cal{G}}}\cap{\cal{G}}^z\cap{\cal{G}}^0\cap{\cal{L}}^{\rho(k)}$ such that (4.1) holds true for ${\cal{R}}$. Since $\widetilde{{\cal{R}}}_T^g({\cal{H}}^{\infty})\in L^2({\cal{F}}_T)$, by Theorem 4.6 again, we further get that $g$ is a function such that (4.1) holds true for $\overline{{\cal{R}}}$. By Lemma A.2(ii) and \cite[Theorem 3.2]{Peng04}, we can get that $L^2$ convergence holds true for ${\cal{E}}[\cdot|{\cal{F}}_t]$ and ${\cal{E}}^g[\cdot|{\cal{F}}_t]$ on $L^2({\cal{F}}_T)$ for each $t\in[0,T]$. Moreover, since the set of all the smooth variables such as
$f(B_{t_1},B_{t_2},\cdots,B_{t_m})$, where $m\geq1$, $f\in C_b^\infty({\bf{R}}^m)$ and $t_i\in[0,T], 1\leq i\leq m,$ is dense in $L^2$ (see \cite[Page 25]{Na}), we conclude that $g$ is a unique function in $\widehat{{\cal{G}}}\cap{\cal{G}}^z\cap{\cal{G}}^0\cap{\cal{L}}^{\rho(k)}$ such that  (4.1) holds true for $L^2({\cal{F}}_T)$.\ \ $\Box$\\\\
\textbf{Remark 4.8} \emph{Theorem 4.6 contains the representation theorems for ${\cal{F}}$-expectations in the locally Lipschitz case, in particular for quadratic ${\cal{F}}$-expectations ($\rho(|\cdot|)=\nu(|\cdot|+1)$). Our study relies on the boundeness of $Z_t$ heavily, which is inspired by the works of \cite{BE} and \cite{CN} depending on Malliavin calculus. Comparing with the representation theorem for quadratic ${\cal{F}}$-expectations in \cite{Yao}, our condition seems to be more simple and not hard to be verified in applications. We can firstly check the dominations on $\cal{R}$ to determine the generator $g$, then further check the dominations on $\widetilde{{\cal{R}}}^{g}({\cal{H}}^{\infty})$. Our proof is more simple, since some subtle techniques of stopping times are omitted by considering a very special ${\cal{E}}$-supermartingale and the set ${\cal{R}}$, which indeed determine the basic part of $g$-expectations. This method also gives a simple proof for the representation theorem in \cite{CH} in the deterministic generator case (see Corollary 4.7).}
\\\\

\begin{center}
\textbf{Appendix: A special nonlinear Doob-Meyer decomposition}
\end{center}

We will introduce a special Doob-Meyer decomposition for ${\cal{E}}$-supermartingales in this section. We firstly recall the definition of ${\cal{E}}$-martingale introduced in \cite{CH}.\\\\
\textbf{Definition A.1} Let ${\cal{E}}$ be an ${\cal{F}}$-expectation. A process $Y_t$ with $Y_t\in L^2({\cal{F}}_t)$ for each $t\in [0,T]$, is called an ${\cal{E}}$-martingale (resp. ${\cal{E}}$-supermartingale, ${\cal{E}}$-submartingale), if for each $0\leq s\leq t\leq T,$ ${\cal{E}}[Y_t|{\cal{F}}_s]=Y_s$ (resp. $\leq,\ \geq$).\\

We recall the domination condition (A$_{\rho(k)}$) introduced in Section 4 and list some properties for the ${\cal{F}}$-expectation ${\cal{E}}$ under (A$_{\rho(k)}$).
\begin{itemize}
  \item (A$_{\rho(k)}$) ($\rho(k)$-domination on ${\cal{K}}$)  For $k>0$, $p>1$ and $\xi, \eta\in{\cal{K}}\subset L^p({\cal{F}}_T)$, we have
\begin{center}
${\cal{E}}[\xi|{\cal{F}}_t]-{\cal{E}}[\eta|{\cal{F}}_t]\leq{\cal{E}}^{\rho(k)}[\xi-\eta|{\cal{F}}_t],\ \ \forall t\in[0,T].$
\end{center}
\end{itemize}
\textbf{Lemma A.2} \emph{Let $p>1$, $k>0$ and the ${\cal{F}}$-expectation ${\cal{E}}$ satisfy (A$_{\rho(k)}$) on ${\cal{K}}\subset L^p({\cal{F}}_T)$. Then for $\xi, \eta\in{\cal{K}}$ and each $t\in[0,T],$ we have}

\emph{(i) ${\cal{E}}^{-\rho(k)}[\xi|{\cal{F}}_t]\leq{\cal{E}}[\xi|{\cal{F}}_t]\leq{\cal{E}}^{\rho(k)}[\xi|{\cal{F}}_t];$}

\emph{(ii) $|{\cal{E}}[\xi|{\cal{F}}_t]-{\cal{E}}[\eta|{\cal{F}}_t]|\leq{\cal{E}}^{\rho(k)}[|\xi-\eta||{\cal{F}}_t];$}

\emph{(iii) ${\cal{E}}[\xi+\zeta|{\cal{F}}_t]={\cal{E}}[\xi|{\cal{F}}_t]+\zeta,$ for $\zeta\in L^p({\cal{F}}_t)$ such that $\xi+\zeta\in{\cal{K}}$.}
\\\\
\emph{Proof.} By (A$_{\rho(k)}$), we have
 $${\cal{E}}^{-\rho(k)}[\xi-\eta|{\cal{F}}_t]=-{\cal{E}}^{\rho(k)}[\eta-\xi|{\cal{F}}_t]\leq{\cal{E}}[\xi|{\cal{F}}_t]-{\cal{E}}[\eta|{\cal{F}}_t]
 \leq{\cal{E}}^{\rho(k)}[\xi-\eta|{\cal{F}}_t].\eqno(a.1)$$
By setting $\eta=0$ in ($a.1$), we get (i). By ($a.1$) and Proposition 3.2(i), we have
$$-{\cal{E}}^{\rho(k)}[|\xi-\eta||{\cal{F}}_t]={\cal{E}}^{-\rho(k)}[-|\xi-\eta||{\cal{F}}_t]
\leq{\cal{E}}[\xi|{\cal{F}}_t]-{\cal{E}}[\eta|{\cal{F}}_t]\leq{\cal{E}}^{\rho(k)}[|\xi-\eta||{\cal{F}}_t],$$
which implies (ii). By ($a.1$) and Proposition 3.3(ii), for each $\zeta\in L^p({\cal{F}}_t)$ such that $\xi+\zeta\in{\cal{K}},$ we have $$\zeta={\cal{E}}^{-\rho(k)}[\zeta|{\cal{F}}_t]\leq{\cal{E}}[\xi+\zeta|{\cal{F}}_t]-{\cal{E}}[\xi|{\cal{F}}_t]\leq{\cal{E}}^{\rho(k)}[\zeta|{\cal{F}}_t]=\zeta.$$
which implies (iii). \ \ $\Box$ \\\\
\textbf{Lemma A.3}\ \ \emph{Let $k>0$ and the ${\cal{F}}$-expectation ${\cal{E}}$ satisfy (A$_{ind}$) and (A$_{\rho(k)}$) on ${\cal{R}}^k$. Then for each $\xi\in{\cal{R}}^k$, ${\cal{E}}[\xi|{\cal{F}}_t]$ is continuous in $t$.}\\\\
\
\emph{Proof.} For each $\xi\in{\cal{R}}^k$, there exist $(y,z)\in{\mathbf{R}}\times{\mathbf{B}}(k)$ and $0\leq u<v\leq T$, such that $\xi=y+z(B_v-B_u).$ By Lemma A.2(iii) and (A$_{ind}$), we have
$${\cal{E}}[\xi|{\cal{F}}_t]={\cal{E}}[y+z(B_v-B_u)|{\cal{F}}_t]=\left\{
                                                                   \begin{array}{ll}
                                                                     y+{\cal{E}}[z(B_v-B_u)],\ \ & t\in[0,u]; \\
                                                                     y+{\cal{E}}[z(B_v-B_t)]+z(B_t-B_u),\ \ & t\in(u,v];\\
                                                                     y+z(B_v-B_u),\ \ & t\in(v,T].
                                                                   \end{array}
                                                                 \right.\eqno(a.2)$$
where implies that for each $t\in[0,T]$, ${\cal{E}}[\xi|{\cal{F}}_t]\in{\cal{R}}^k.$ Then by ($a.2$), Lemma A.2(ii) and \cite[Theorem 3.2]{Peng04}, we deduce that ${\cal{E}}[\xi|{\cal{F}}_t]$ is continuous in $t$.\ \  $\Box$ \\\\
\textbf{Lemma A.4} \emph{Let $k>0$ and the ${\cal{F}}$-expectation ${\cal{E}}$ satisfy (A$_{ind}$) and (A$_{\rho(k)}$) on ${\cal{R}}^k$. Then for each $\xi\in{\cal{R}}^k$, there exists a pair $(g_t^\xi, Z_t^\xi)\in{\cal{H}}^2_{1\times 1}\times{\cal{H}}^2_{1\times d}$ such that} $$|g_t^\xi|\leq \rho(k)|Z_t^\xi|,\ dt\times dP-a.e., \ \emph{and}\ \
{\cal{E}}[\xi|{\cal{F}}_t]=\xi+\int_t^Tg_s^\xi ds-\int_t^TZ_s^\xi dB_s,\ \ t\in[0,T].\eqno(a.3)$$
\emph{Moreover, for $\eta\in{\cal{R}}^k$, we have} $$|g_t^\xi-g_t^\eta|\leq\rho(k)|Z_t^\xi-Z_t^\eta|,\ \ dt\times dP-a.e.\eqno(a.4)$$
\emph{Proof.} From Lemma A.3 and the proof of \cite[Lemma 5.3]{CH}, we can get ($a.3$) directly. Now we prove ($a.4$).
By ($a.3$), for $\eta\in{\cal{R}}^k$, there exists a pair $(g_t^\eta, Z_t^\eta)\in{\cal{H}}^2_{1\times 1}\times{\cal{H}}^2_{1\times d}$ such that
$${\cal{E}}[\eta|{\cal{F}}_t]=\eta+\int_t^Tg_s^\eta ds-\int_t^TZ_s^\eta dB_s,\ \ t\in[0,T].\eqno(a.5)$$
Set $y_t^\xi:={\cal{E}}[\xi|{\cal{F}}_t]$ and $y_t^\eta:={\cal{E}}[\eta|{\cal{F}}_t]$. By ($a.2$), we have $y_t^\xi, y_t^\eta\in{\cal{R}}^k$ for each $t\in[0,T]$. Then by ($a.3$), ($a.5$) and ($a.1$), we have for each $0\leq s\leq t\leq T,$
$$y_s^\xi-y_s^\eta={\cal{E}}[y_t^\xi|{\cal{F}}_s]-{\cal{E}}[y_t^\eta|{\cal{F}}_s]\geq{\cal{E}}^{-\rho(k)}[y_t^\xi-y_t^\eta|{\cal{F}}_s]$$
and
$$y_s^\xi-y_s^\eta={\cal{E}}[y_t^\xi|{\cal{F}}_s]-{\cal{E}}[y_t^\eta|{\cal{F}}_s]\leq{\cal{E}}^{\rho(k)}[y_t^\xi-y_t^\eta|{\cal{F}}_s].$$
Thus, $y_t^\xi-y_t^\eta$ is a $\rho(k)$-submaringale and a $-\rho(k)$-supermaringale. Then by a similar argument as in ($a.3$) (see \cite[Lemma 5.3]{CH}), we can get ($a.4$). \ \  $\Box$\\

Let the function $f:[0,T]\times \mathbf{R}\longmapsto \mathbf{R}$ such that $f(\cdot,0)\in l^2_{1\times1}(0,T)$ and for each $y_1, y_2\in{\mathbf{R}}$,
$$|f(t,y_1)-f(t,y_2)|\leq \lambda|y_1-y_2|,\ \ \forall t\in[0,T],$$
where $\lambda>0$ is a constant. For $(y,z)\in\textbf{R}\times\textbf{B}(k)$, we consider the following BSDE ${\cal{E}}(f,T,y,z):$
$$\left\{
      \begin{array}{ll}
        \psi(t)+zB_t={\cal{E}}\left[y+zB_T+\int_t^Tf(s,\psi(s))ds|{\cal{F}}_t\right], \ \ t\in[0,T); \\
        \psi(T)=y,\ \ y\in \textbf{R}.
      \end{array}
    \right.$$
\textbf{Lemma A.5} \emph{Let $k>0$ and the ${\cal{F}}$-expectation ${\cal{E}}$ satisfy (A$_{ind}$) and (A$_{\rho(k)}$) on ${\cal{R}}^k$. Then for each $(y,z)\in{\bf{R}}\times{\bf{B}}(k)$, the BSDE ${\cal{E}}(f,T,y,z)$ has a unique solution in $l^2_{1\times1}(0,T).$} \\\\
\emph{Proof.} For $\psi(t)\in l^2_{1\times1}(0,T),$ set $$I(\psi(t)):={\cal{E}}\left[y+zB_T+\int_t^Tf(s,\psi(s))ds|{\cal{F}}_t\right]-zB_t,\ \ t\in[0,T].$$
Since $f(\cdot,0)\in l^2_{1\times1}(0,T),$ we have $$\int_0^T|f(s,\psi(s))|^2ds\leq2\int_0^T(|f(s,0)|^2+\lambda^2|\psi(s)|^2)ds<\infty.$$
Then by Lemma A.2(iii), (A$_{ind}$), Lemma A.2(i) and (2.17), we can get that
\begin{eqnarray*}
|I(\psi(t))|&=&|{\cal{E}}[z(B_T-B_t)]|+|y|+|\int_t^Tf(s,\psi(s))ds|\\
&\leq&|{\cal{E}}^{-\rho(k)}[z(B_T-B_t)]|+|{\cal{E}}^{\rho(k)}[z(B_T-B_t)]|+|y|+\int_t^T|f(s,\psi(s))|ds\\
&\leq&2\rho(k)|z|T+|y|+\int_t^T|f(s,\psi(s))|ds,\ \ t\in[0,T],
\end{eqnarray*}
which implies that $I(\psi(t))\in l^2_{1\times1}(0,T)$, i.e. $$I(\cdot):l^2_{1\times1}(0,T)\longmapsto l^2_{1\times1}(0,T).$$

Since for each $\psi(t)\in l^2_{1\times1}(0,T),$ $y+zB_T+\int_t^Tf(s,\psi(s))ds\in{\cal{R}}^k$. Then by Lemma A.2(ii) and Proposition 3.3(ii), for $\psi^1(t), \psi^2(t)\in l^2_{1\times1}(0,T),$ we have
\begin{eqnarray*}
  &&|I(\psi^1(t))-I(\psi^2(t))|^2\\&=&\left|{\cal{E}}\left[y+zB_T+\int_t^Tf(s,\psi^1(s))ds|{\cal{F}}_t\right]
  -{\cal{E}}\left[y+zB_T+\int_t^Tf(s,\psi^2(s))ds|{\cal{F}}_t\right]\right|^2\\
  &\leq&\left|{\cal{E}}^{\rho(k)}\left[\left|\int_t^Tf(s,\psi^1(s))ds
  -\int_t^Tf(s,\psi^2(s))ds\right||{\cal{F}}_t\right]\right|^2\\
  &\leq&\left|\int_t^T\left|f(s,\psi^1(s))-f(s,\psi^2(s))\right|ds\right|^2\\
  &\leq&\lambda^2T^2 \int_t^T|\psi^1(s)-\psi^2(s)|^2ds,\ \ t\in[0,T].
\end{eqnarray*}
It follows that
\begin{eqnarray*}
\int_0^T|I(\psi^1(t))-I(\psi^2(t))|^2e^{2\lambda^2T^2t}dt&\leq&\lambda^2T^2\int_0^Te^{2\lambda^2T^2t}\int_t^T|\psi^1(s)-\psi^2(s)|^2dsdt\\
&=&\lambda^2T^2\int_0^T\int_0^se^{2\lambda^2T^2t}|\psi^1(s)-\psi^2(s)|^2dtds\\
&=&\frac{1}{2}\int_0^T(e^{2\lambda^2T^2s}-1)|\psi^1(s)-\psi^2(s)|^2ds\\
&\leq&\frac{1}{2}\int_0^T|\psi^1(s)-\psi^2(s)|^2e^{2\lambda^2T^2s}ds,
\end{eqnarray*}
which means that $I(\cdot)$ is a strict contraction on $l^2_{1\times1}(0,T)$ under the norm $\int_0^T|\cdot|e^{2\lambda^2T^2s}ds$. The proof is complete.\ \  $\Box$\\

By Lemma A.3, for each $(y,z)\in{\bf{R}}\times{\bf{B}}(k)$, the solution $\psi(t)$ of the BSDE ${\cal{E}}(f,T,y,z)$ is continuous. Now, we give a special Doob-Meyer decomposition of ${\cal{E}}$-supermartingale.\\\\
\textbf{Lemma A.6} \emph{Let $k>0$ and the ${\cal{F}}$-expectation ${\cal{E}}$ satisfy (A$_{ind}$) and (A$_{\rho(k)}$) on ${\cal{R}}^k$. If $\psi(t)\in l^2_{1\times1}(0,T)$ and $z\in{\mathbf{B}}(k)$ such that $\psi(t)+zB_t$ is a continuous ${\cal{E}}$-supermartingale, then there exists a function $a(t)\in l^2_{1\times1}(0,T)$, which is continuous and increasing with $a(0)=0$ such that}
$${\cal{E}}[\psi(T)+zB_T+a(T)|{\cal{F}}_t]=\psi(t)+zB_t+a(t),\ \ \forall t\in[0,T].$$
\emph{Proof.} By Lemma A.5, for each $m\geq1$, the following equation
$$\psi^m(t)+zB_t={\cal{E}}\left[\psi(T)+zB_T+\int_t^Tm(\psi(s)-\psi^m(s))ds|{\cal{F}}_t\right],\ \ t\in[0,T],\eqno(a.6)$$
has a unique solution $\psi^m(t)\in l^2_{1\times1}(0,T)$. We will prove that for each $m\geq1$ and $t\in[0,T]$, $\psi(t)\geq\psi^m(t)$. For $m\geq1$, if there exists $r\in[0,T]$ such that $\psi(r)<\psi^m(r)$, then by the continuity of $\psi(t)$ and $\psi^m(t)$, and the fact $\psi(T)=\psi^m(T)$, there exists $s\in(r,T]$ such that for each $t\in[r,s)$, $\psi(t)<\psi^m(t)$ and $\psi(s)=\psi^m(s).$ Since $\psi(t)+zB_t$ is an ${\cal{E}}$-supermartingale, then by Lemma A.2(iii), for $t\in[r,s)$, we have
\begin{eqnarray*}
\psi(t)+zB_t&\geq&{\cal{E}}\left[\psi(s)+zB_s|{\cal{F}}_t\right]\\&=&{\cal{E}}\left[\psi^m(s)+zB_s|{\cal{F}}_t\right]\\
&>&{\cal{E}}\left[\psi^m(s)+zB_s+\int_t^sm(\psi(u)-\psi^m(u))du|{\cal{F}}_t\right]\\
&=&{\cal{E}}\left[{\cal{E}}\left[\psi(T)+zB_T+\int_s^Tm(\psi(u)-\psi^m(u))du|{\cal{F}}_s\right]+\int_t^sm(\psi(u)-\psi^m(u))du|{\cal{F}}_t\right]\\
&=&{\cal{E}}\left[\psi(T)+zB_T+\int_t^Tm(\psi(u)-\psi^m(u))du|{\cal{F}}_t\right]\\
&=&\psi^m(t)+zB_t
\end{eqnarray*}
which implies $\psi(t)>\psi^m(t)$. This leads to a contradiction. Consequently, for each $m\geq1$ and $t\in[0,T]$, we have $$\psi(t)\geq\psi^m(t).$$ By Lemma A.2(iii) again, we have
\begin{eqnarray*}
\psi^{m+1}(t)-\psi^m(t)&=&{\cal{E}}\left[\psi(T)+zB_T+\int_t^T(m+1)(\psi(u)-\psi^{m+1}(u))du|{\cal{F}}_t\right]\\
\ \ \ \ \ \ &&-{\cal{E}}\left[\psi(T)+zB_T+\int_t^Tm(\psi(u)-\psi^m(u))du|{\cal{F}}_t\right]\\
&=&-\int_t^Tm(\psi^{m+1}(u)-\psi^m(u))du+\int_t^T(\psi(u)-\psi^{m+1}(u))du, \ \ \ t\in[0,T].
\end{eqnarray*}
Since $\psi(u)-\psi^{m+1}(u)\geq0$, by solving the above linear ODE, we can get that for each $t\in[0,T]$, $\psi^{m+1}(t)\geq\psi^m(t).$ By setting
$$a^m(t):=\int_0^tm(\psi(s)-\psi^m(s))ds, \ \ \ t\in[0,T],\eqno(a.7)$$
and ($a.6$), we have $$\psi^m(t)={\cal{E}}\left[\psi(T)+zB_T|{\cal{F}}_t\right]-zB_t+a^m(T)-a^m(t), \ \ \ t\in[0,T].\eqno(a.8)$$
Since as $m\rightarrow\infty,$ for each $t\in[0,T]$, the limit of $\psi^m(t)$ exists, it follows from ($a.8$) that for each $t\in[0,T]$, the limit of $a^m(T)-a^m(t)$ exists, which together with the fact $a^m(0)=0$ implies that for each $t\in[0,T]$, the limit of $a^m(t)$ exists. We denote the limit of $a^m(t)$ by $\tilde{a}(t).$ Then by ($a.7$), we can get that for almost each $t\in[0,T]$, $\psi^m(t)\rightarrow\psi(t)$, as $m\rightarrow\infty,$ and $\tilde{a}(t)$ is a deterministic and increasing function with $\tilde{a}(0)=0.$ Then by letting $m\rightarrow\infty$ in ($a.8$), we have
$$\psi(t)={\cal{E}}\left[\psi(T)+zB_T|{\cal{F}}_t\right]-zB_t+\tilde{a}(T)-\tilde{a}(t), \ \ \ a.e.\ t\in[0,T].$$
Moreover, since $\psi(t)$ is continuous and $\tilde{a}(t)$ is monotonic, by setting $\hat{a}(t):=\lim_{s>t, s\rightarrow t}\tilde{a}(s)$, $t\in[0,T)$, and $\hat{a}(T):=\tilde{a}(T)$, we have
$$\psi(t)={\cal{E}}[\psi(T)+zB_T|{\cal{F}}_t]-zB_t+a(T)-a(t),\ \ \textmd{with} \ a(t):=\hat{a}(t)-\hat{a}(0),\  \forall t\in[0,T].\eqno(a.9)$$
It is clear that $a(t)\in l^2_{1\times1}(0,T)$ is deterministic, continuous and increasing with $a(0)=0$. Then by Lemma A.2(iii) and ($a.9$), we can complete this proof. \ \ $\Box$\\\\
\textbf{Declaration of competing interest}\\

The authors declare that they have no known competing financial interests or personal
relationships that could have appeared to influence the work reported in this paper.

\end{document}